\documentclass{article}
\usepackage{amssymb,amsthm,amsmath}
\usepackage[hmargin=2.5cm,vmargin=3.8cm]{geometry}
\usepackage[T1]{fontenc}
\usepackage[utf8]{inputenc}         

\usepackage{float}

\usepackage{graphicx} 
\usepackage{enumerate}
\usepackage{enumitem, hyperref}
\usepackage{cleveref}
\usepackage{xspace}
\usepackage{subcaption}

\usepackage{authblk}
\usepackage{todonotes}
\usepackage{comment}

\usepackage{tikz}
\tikzstyle{v}=[circle, inner sep=1pt, minimum size =18 pt, line width = 1pt, draw=black, fill=white, text= black]

\usepackage{algorithm,algpseudocode}

\theoremstyle{plain}
\newtheorem{theorem}{Theorem}

\newtheorem{lemma}[theorem]{Lemma}
\newtheorem{corollary}[theorem]{Corollary}
\theoremstyle{definition}
\newtheorem{definition}[theorem]{Definition}
\newtheorem{example}[theorem]{Example}
\newtheorem{remark}[theorem]{Remark}

\newtheorem{claim}{Claim}[theorem]
\newtheorem*{notation}{Notation}

\Crefname{claim}{Claim}{Claims}

\newcommand{\G}{\mathcal{G}}
\newcommand{\M}{\mathcal{M}}
\newcommand{\Pc}{\mathcal{P}}
\DeclareRobustCommand{\rchi}{{\mathpalette\irchi\relax}}
\newcommand{\irchi}[2]{\raisebox{\depth}{$#1\chi$}} 

\title{Reconstructing graphs with subgraph compositions\footnote{The research of the second author was partially funded by the Academy of Finland grants 338797 and 358718.}~\footnote{A shorter version of this paper has been submitted to a conference.}}
\author[1]{Antoine Dailly}
\author[2]{Tuomo Lehtilä}
\date{}
\affil[1]{Université Clermont Auvergne, INRAE, UR TSCF, 63000, Clermont-Ferrand, France}
\affil[2]{Department of Mathematics and Statistics, University of Turku, FI-20014, Turku, Finland}

\begin{document}

\maketitle

\begin{abstract}
	We generalize the problem of reconstructing strings from their substring compositions first introduced by Acharya \emph{et al.} in 2015 motivated by polymer-based advanced data storage systems utilizing mass spectrometry. Namely, we see strings as labeled path graphs, and as such try to reconstruct labeled graphs. For a given integer $t$, the subgraph compositions contain either vectors of labels for each connected subgraph of order~$t$ ($t$-multiset-compositions) or the sum of all labels of all connected subgraphs of order~$t$ ($t$-sum-composition). We ask whether, given a graph of which we know the structure and an oracle whom you can query for compositions, we can reconstruct the labeling of the graph. If it is possible, then the graph is reconstructable; otherwise, it is confusable, and two labeled graphs with the same compositions are called equicomposable. We prove that reconstructing through a brute-force algorithm is wildly inefficient, before giving methods for reconstructing several graph classes using as few compositions as possible. We also give negative results, finding the smallest confusable graphs and trees, as well as families with a large number of equicomposable non-isomorphic graphs. An interesting result occurs when twinning one leaf of a path: some paths are confusable, creating a twin out of a leaf sees the graph alternating between reconstructable and confusable depending on the parity of the path, and creating a false twin out of a leaf makes the graph reconstructable using only sum-compositions in all cases.
\end{abstract}


\section{Introduction}

The problem of reconstructing a graph from incomplete information has been studied for quite some time, being introduced at least in the 1950s~\cite{kelly1957congruence}. Much of the research on this topic is about finding the original graph by asking an oracle questions about, \emph{e.g.}, connectivity of sets of vertices. In particular, we usually know the vertex set and wish to find the edge-set. 
We are interested in a vertex-\emph{labeled} version of the reconstruction problem: given a graph (of which we know the structure) and an oracle, can we find the labels of the vertices? Note that we are in a wholly deterministic setting, so the oracle answers truthfully and correctly, and our goal is to find the correct labels exactly (up to graph isomorphisms).
Furthermore, the problem is trivial if only one label can be used, so we are going to use alphabets of size at least~2.

There are several ways to define an oracle to query about the graph, and we are going to use the one used by Acharya \emph{et al.} in the context of strings~\cite{acharya2015string} and reused by Bartha \emph{et al.} for edge-labeled graphs~\cite{bartha2021reconstructibility,bartha2016reconstruction}. The oracle can be queried about \emph{subgraph compositions}: for a given integer $t$, it gives information about all connected subgraphs of order $t$. This information is not separated, so we do not know which piece of information relates to any given subgraph.

Given a vertex-labeled graph $\G$ where each label belongs to $k$-letter alphabet $\Sigma_k$ for $k\geq2$, we give the total number each label occurs with an integer valued vector of length $k$. For example, if the graph has two vertices labeled with $A$ and one with $B$ over alphabet $\Sigma_2=(A,B)$, then our vector would be $(2,1)$.
We are going to use two types of compositions. First is the \emph{multiset-composition} $\mathcal{M}_t(\mathcal{G})$, which is a multiset of vectors such that each vector corresponds to a connected subgraph of order $t$, and contains all of its labels (note that we do not know which vector is connected to which subgraph, so we do not know the labels of vertices but only how many labels there are in different types of subgraphs). Next, the \emph{sum-composition} $S_t(\mathcal{G})$, which is a vector of length $k$, where each element corresponds to a symbol of the alphabet $\Sigma_k$, and contains the number of times the symbol appears in all subgraphs of order $t$. It is clear that we can find a sum-composition from a multiset-composition, but the reverse is not true, hence sum-compositions contain less information, and being able to reconstruct a graph using only sum-compositions is generally stronger.

If it is possible to connect all labels of a graph $G$ with its vertices by querying compositions, then we say that $G$ is \emph{reconstructable}. There is a straightforward and wildly inefficient algorithm to find the labels of a reconstructable graph: enumerate all its multiset-compositions, and do so for every possible non-isomorphic labeling until the multiset-combinations coincide. Obviously, our goal is to have more efficient algorithms for reconstructing graphs, which we will generally evaluate by number of queries (note that, in a classical complexity framework, we consider a query as a singular operation, having access to the oracle). Since sum-compositions contain less information than multiset-compositions, we also consider algorithms reconstructing a graph using only sum-compositions to be better than those also using multiset-compositions. If $G$ can be reconstructed with only sum-compositions, then we say that $G$ is \emph{sum-reconstructable}.

Some graphs cannot be reconstructed, even with the brute-force method, those are called \emph{confusable}. Two labeled graphs that yield the exact same multiset-compositions are called \emph{equicomposable}, and a simple way to prove that a graph $G$ is confusable is to exhibit two equicomposable labeled copies of $G$.

Formal definitions will be given in \Cref{subsec-definitions}. Before this, we explain in the following subsection how the problem of reconstructing graphs can be related to recovering data from polymers.

\subsection{Polymer-based memory systems}\label{SubsecPolymer}

Recently, polymer-based data storage systems have been considered from both theoretical and experimental points-of-views, see for example~\cite{acharya2015string, al2017mass, launay2021precise}. In these types of memory systems, the information is stored into 
the (synthesized) polymer using molecules with large weight differences.  Then, for example, when we use molecules with weights $A$ and $B$, we can deduce the number of $A$'s and $B$'s based on the total weight of the polymer but this does not allow us to immediately know their ordering. For example, when we know that there are at most nine molecules in total and we have $A=1$ and $B=10$, then total mass 72 implies that there are two $A$s and seven $B$s.

For reading the mass of the polymer,  tandem mass (MS/MS) spectrometers are used \cite{acharya2015string,gabrys2022reconstruction, launay2021precise}. The tandem mass spectrometers break the polymers into smaller fragments and measure the masses of these fragments. The goal is to reconstruct the original polymer using the masses of these fragments which we obtain from multiple original copies of the polymer.

We consider this problem from a somewhat idealized combinatorial perspective for mathematical modeling, as has been previously done for example in \cite{acharya2015string,gabrys2020mass,gabrys2022reconstruction, pattabiraman2019reconstruction,  pattabiraman2023coding}. The two assumptions originally presented in~\cite{acharya2015string} which we use are:
\begin{itemize}
	\item[A1:] We can uniquely deduce the composition, that is, the number of each symbol of a fragment string from its mass.
	\item[A2:] We observe every fragment of a polymer with equal frequency.
\end{itemize}
In previous work, it has been assumed that the underlying polymer structure is a string (although in \cite{banerjee2023insertion} generalizing it to other structures, such as cycles, has been suggested).

\subsection{Problem definitions and first lemmas}
\label{subsec-definitions}

In this section, we give formal definitions, explain our notation and present some general lemmas.

\begin{definition}\label{defMultiset}
	For an integer $k \geq 2$, let $\mathcal{G}=(G,\lambda)$ be a connected labeled graph with labeling $\lambda:V(G)\to \Sigma_k$, where $\Sigma_k$  contains symbols $\{A_1,\dots,A_k\}$. We denote by $n_i(\mathcal{G})=|\{v\in V(G)\mid \lambda(v)=A_i\}|$ and by $c(\mathcal{G})=\{A_1^{n_1(\mathcal{G})},A_2^{n_2(\mathcal{G})},\dots,A_k^{n_k(\mathcal{G})}\}$.
	
	We denote the \emph{subgraph composition multiset} of $\mathcal{G}$ by: $$\mathcal{M}(\mathcal{G})=\{\{c(\mathcal{G'})\mid \mathcal{G'}\text{ is connected induced subgraph of } \mathcal{G} \}\}$$ and call the elements of $\mathcal{M}(\mathcal{G})$ \emph{compositions}.

	Let $\mathcal{G}_t$ be the set of $t$-vertex connected induced subgraphs of $\mathcal{G}$. We denote by $\M_t(\G)$ the restriction of $\M(\G)$ to connected $t$-vertex subgraphs of $\G$, that is, $\M_t(\G)=\{\{c(\mathcal{G'})\mid \mathcal{G'}\in \mathcal{G}_t \}\}$. We define the \emph{$t$-sum-composition} of $\mathcal{G}$ as $S_t(\mathcal{G})=\sum_{\mathcal{G}'\in \mathcal{G}_t}(n_1(\mathcal{G}'),n_2(\mathcal{G}'),\dots,n_k(\mathcal{G}'))$ and the sum-composition of $\mathcal{G}$ by $S(\mathcal{G})=\bigcup_{t=1}^{|V(G)|} S_t(\mathcal{G})$.
\end{definition}

\begin{notation}
	For readability, when $c(\mathcal{G})=\{A_1^{n_1(\mathcal{G})},A_2^{n_2(\mathcal{G})},\dots,A_k^{n_k(\mathcal{G})}\}$, we may denote it as a sequence, that is, $c(\mathcal{G})=A_1^{n_1(\mathcal{G})}A_2^{n_2(\mathcal{G})} \ldots A_k^{n_k(\mathcal{G})}$.
\end{notation}

\begin{notation}
	For ease of use, we will sometimes denote the $t$-sum-composition of $\mathcal{G}$ as:
	\[ S_t(\mathcal{G}) = \sum_{i=1}^k n_i^t(\mathcal{G}) A_i \]
	where $n_i^t(\mathcal{G}) = \sum_{\mathcal{G}' \in \mathcal{G}_t} n_i(\mathcal{G}')$. That is, $n_i^t(\mathcal{G})$ is the number of times the label $A_i$ appears in total in every connected induced subgraph of order $t$.
\end{notation}

Note that $c(\mathcal{G})$ can be understood as a \textit{Parikh vector} of $\mathcal{G}$,~\cite{bartha2021reconstructibility}. Hence we mean with, for example, alphabet $\Sigma_3=(A,B,C)$ that composition $A$ is the same as vector $(1,0,0)$ and $(2,1,0)-AB=(1,0,0)=A$ (the subtraction $X-Y$ is only allowed if it does not result in any negative coordinates).
We consider $t$-subgraph composition multisets and $t$-sum compositions as \emph{queries} asked to an oracle. We call the former a \textit{multiset-query} and the latter a \textit{sum-query}.

\begin{figure}[h]
	\centering
	\begin{tikzpicture}
		\node[v] (c) at (3,1.5) {$B$};
		\node[v] (0) at (1.5,1.5) {$A$};
		\node[v] (1) at (0,0) {$A$};
		\node[v] (2) at (0,3) {$C$};
		\node[v] (3) at (4.5,1.5) {$D$};
		\node[v] (4) at (6,0) {$A$};
		\node[v] (5) at (6,3) {$D$};
		\node[v] (6) at (3,0) {$C$};
		\node[v] (7) at (3,3) {$D$};
		\draw (4)to(3)to(c)to(0)to(1);
		\draw (0)to(2);
		\draw (3)to(5);
		\draw (6)to(c)to(7);
	\end{tikzpicture}
	\caption{A labeled tree $\mathcal{T}$ on nine vertices.}
	\label{fig:Ex1graph}
\end{figure}
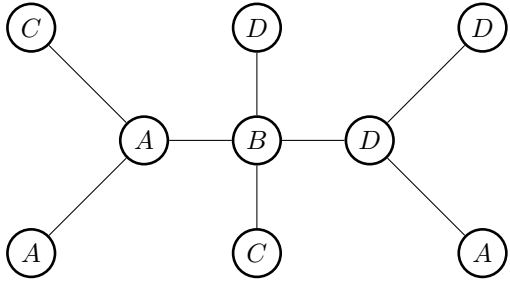

As an example, consider the labeled tree $\mathcal{T}$ from \Cref{fig:Ex1graph}. The multiset-query $\mathcal{M}_3(\mathcal{T})$ outputs $\{\{ A^2C,$ $ABC,$ $ ABC, A^2B, ABD, ABD, ABD, BD^2, BD^2, BCD, BCD, AD^2 \}\}=\{\{ A^2C,$ $2ABC, A^2B, 3ABD,$ $ 2BD^2,$ $ 2BCD, AD^2 \}\}$ where we write $A^2 $ for $AA$ and $2ABC$ for two copies of $ABC$. Note that in our notation $ABC$ and $ACB$ are the same as we are only interested in the compositions. The sum-query $S_3(\mathcal{T})$ outputs $(10,10,5,11)$, which can also be written as $S_3(\mathcal{T})=10A+10B+5C+11D$.

In the following definition, we define what we mean by \textit{reconstructable}. Informally, the question is: can we (or \emph{how} can we) obtain the original graph labeling of a graph (up to graph isomorphisms), when we know the structure of the underlying graph and the multiset $\M(\G)$. Furthermore, we also define a stronger notion of \textit{sum-reconstructability} which asks whether the labeling of $\G$ can be deduced based on $S(\mathcal{\G})$ instead of $\M(\G)$.

\begin{definition}\label{defReconstructable}
	Given a graph $G$, let $\G_1=(G,\lambda_1)$ and $\mathcal{G}_2=(G,\lambda_2)$ be two non-isomorphic labeled graphs.
	
	We say that $G$ is \emph{reconstructable} if we have $\mathcal{M}(\mathcal{G}_1)\neq \mathcal{M}(\mathcal{G}_2)$ for any two non-isomorphic labelings. Similarly, a labeled graph $\mathcal{G}=(G,\lambda)$ is reconstructable if for each non-isomorphic labeled graph $\mathcal{G}'=(G,\lambda')$ we have  $\mathcal{M}(\mathcal{G})\neq \mathcal{M}(\mathcal{G}')$. 
	We say that $G$ is \emph{sum-reconstructable} if we have $S(\mathcal{G}_1)\neq S(\mathcal{G}_2)$.
	
	If $G$ or $\mathcal{G}$ is not reconstructable, then we say that it is \emph{confusable}.
	Furthermore, we say that $\mathcal{G}_1$ and $\mathcal{G}_2$ are \emph{equicomposable} if $\mathcal{M}(\mathcal{G}_1)= \mathcal{M}(\mathcal{G}_2)$.
\end{definition}

As a first, easy result, consider the complete graph $K_n$ on $n$ vertices, where every possible edge exists. Since all vertices are symmetrically positioned, two labelings with the same subgraph composition multiset are isomorphic, and thus it suffices to query $\mathcal{M}_n(K_n)$ to reconstruct $K_n$.

We also show the following natural lemmas:

\begin{lemma}\label{lemSumRecToRec}
	If $G$ is sum-reconstructable, then $G$ is reconstructable.
\end{lemma}

\begin{proof}
	Let $G$ be sum-reconstructable and suppose to the contrary that $\mathcal{G}=(G,\lambda)$ and $\mathcal{G}'=(G,\lambda')$ are equicomposable. Hence, $\M(\G)=\M(\G')$ and in particular, the subset of $\M(\G)$ containing exactly the $t$ symbol compositions of $\M(\G)$ is equal to  the subset of $\M(\G')$ containing exactly the $t$ symbol compositions of $\M(\G')$. Hence, we have $S_t(\G)=S_t(\G')$. Furthermore, this is true for any $t\geq1$. Hence, $S(\G)=S(\G')$, a contradiction. Thus, graph $G$ is reconstructable. 
\end{proof}

\begin{lemma}
	\label{lem-s2G}
	For $\mathcal{G}=(G,\lambda)$, we have $S_2(\mathcal{G})=\sum_{v \in V(G)} d(v) \lambda(v)$.
\end{lemma}

\begin{proof}
	In $S_2(\G)$, every edge is counted. Hence, the label of each vertex is counted exactly in as many edges as the vertex appears, which is its degree.
\end{proof}

\subsection{Graph notations}

The graphs we consider are all simple, connected, and unweighted. We are going to study specific graph families, focusing mainly on \emph{trees}, that is, acyclic graphs. In a tree, degree~1 vertices are called \emph{leaves} while the other vertices are called \emph{internal vertices}. The simplest subfamily of trees is \emph{paths}: the path $P_n$ on $n$ vertices has two leaves (or one if $n=1$). Another simple family is \emph{stars}: the star $S_n$ on $n$ vertices has $n-1$ leaves connected to the \emph{center}. The \emph{subdivided star} $S_{n_1,\ldots,n_m}$ is the star $S_m$ where each edge has been subdivided $n_1-1,\ldots,n_m-1$ times, that is, it consists of a center to which are attached paths of $n_1,\ldots,n_m$ vertices, called \emph{branches}. Finally, a \emph{bistar} $B_{m,n}$ has two adjacent vertices each being adjacent with $m$ and $n$ leaves, respectively.

The \emph{open neighbourhood} of a vertex is the set of all the vertices it is adjacent to, and its \emph{closed neighbourhood} is its open neighbourhood plus itself; we call \emph{twins} (resp. \emph{false twins}) two vertices with the same closed (resp. open) neighbourhood (for instance, in a star, all leaves are false twins, and in the complete graph all vertices are twins). Note that, if a set of vertices are all twins or false twins, then they are not distinguishable in the setting we are studying, hence, we consider isomorphic labelings to be the same. We denote by $\sqcup$ the additive union of (multi)sets, that is, $A\sqcup A= \{\{A,A\}\}$.

\subsection{Related research}\label{subsecRelated}

The related research consists of two topics: graph theoretical reconstruction problems and mass spectrometry related string reconstruction problems as well as DNA-based data storages. We first discuss about  graph theoretic reconstruction problems and then continue with  mass spectrometry and DNA related research.

Probably the best known and still open graph reconstruction conjecture from Kelly and Ulam states that any (non-labeled) graph $G$ on at least $n\geq3$ vertices can be reconstructed (up to isomorphisms) based on the multiset containing every induced $(n-1)$-vertex subgraph  of $G$~\cite{kelly1957congruence, ulam1960collection}. While this conjecture is still open, it has been  proved for some graph classes such as trees~\cite{kelly1957congruence}. Interested readers may read more about the reconstruction conjecture in a survey article by Harary~\cite{harary2006survey}. For some variants of graph reconstruction, we refer interested readers to the survey~\cite{asciak2010survey} or to~\cite{brown2018new, kostochka2020reconstruction, spinoza2019reconstruction}.

Bartha \emph{et al.} have considered an edge version of our reconstruction problem in~\cite{bartha2021reconstructibility, bartha2016reconstruction}. In their variant, the edges are labeled (instead of vertices) and we obtain a composition multiset containing the information about the edge labels. However, their results on the edge-variant do not seem to give us much information about the vertex-labeled case, although they do note some minor results for vertex-labeled case in their discussions. In particular, they state (without a proof or a construction) that the smallest confusable vertex-labeled tree has seven vertices. We obtained the mentioned graph through personal communications with Burcsi and Lipták (see Figure~\ref{fig-s123}) and show its unicity (see \Cref{thm-smallestConfusableGraphs}).

\medskip

From the application side, the research related to advanced polymer-based memory storages is interesting. While especially DNA-based memory storage systems are widely studied, our research is related more strongly to polymer-based data storage systems using mass spectrometry. This has been also the  motivation of Acharya et al. (and the related work it has generated) in~\cite{acharya2015string} who initiated the study of string reconstruction using substring compositions, which we generalize to graphs: while they construct strings, it is possible to interpret their work as research on graph reconstruction focusing on labeled paths.
Real-world mass spectrometry experiments on polymer-based data storage have been considered, for example in~\cite{al2017mass, launay2021precise}.

In particular, Acharya et al. have shown that the path $P_n$ is reconstructable (from $\mathcal{M}(\mathcal{P}_n)$) if and only if $n+1\leq7$ or $n+1\in\{p,2p\}$ where $p$ is a prime. The case of paths is natural from the point of view of memories, as information is often considered as strings of symbols. Banerjee et al.~\cite{banerjee2023insertion} have suggested generalizing the problem of string reconstruction to more complicated structures. In particular, they suggested cycles. 
The research of Acharya et al. has sprung up multiple studies and variations,  see for example~\cite{acharya2015string,  gabrys2020mass, gabrys2022reconstruction, pattabiraman2019reconstruction, pattabiraman2023coding, yang2024reconstruction, ye2023reconstruction}. 

Another related widely studied polymer-based advanced information storage system is DNA-based data storage system where the information is directly stored into DNA-strands. These systems have also been studied from a similar combinatorial reconstruction perspective, where a strand is reconstructed from multiple erroneous copies, see for example~\cite{junnila2021levenshtein, junnila2024levenshtein, Levenshtein, levenshtein2001efficient, PHAM2025105980, YBiram}.
It has been pointed out in~\cite{al2017mass} that the advantage of polymer-based systems over the DNA-based is their lower cost. 

Some of variants of the reconstruction problem include reconstructing strings/paths from composition multisets containing one of the two leaves of the path, that is, reconstructions of strings from prefixes and suffixes~\cite{gabrys2022reconstruction, yang2024reconstruction, ye2023reconstruction} (however, this variant seems most natural in the special case of paths). Another variant suggested by Bartha et al. in~\cite{bartha2021reconstructibility, bartha2016reconstruction} for the edge-labeled reconstruction is to consider instead of any induced subgraph, only subgraphs with specific structure. For example, only subgraphs which are paths. Finally, one more variant, considered by Bartha et al. in~\cite{bartha2021reconstructibility, bartha2016reconstruction} for the edge-labeled case, is to not give the structure of the underlying graph (or to restrict it to some graph class such as trees) and to also reconstruct it.

\subsection{Structure of the paper}

We begin in \Cref{SecLabels} by showing that the brute-force reconstruction can be extremely inefficient, since the graphs we are studying have a large number of non-isomorphic labelings. This also relates to the application to information storage. Then, in \Cref{sec-reconstructing}, we show how to reconstruct some classes of trees. An interesting case arises when creating a twin from a leaf of a path: this simple operation has a significant impact on reconstructability. The choice of creating a twin or a false twin also has an impact, as is shown in \Cref{secSmallConfusable}, where we also enumerate the smallest confusable graphs and trees. Finally, in \Cref{secConfusable}, we construct large families of equicomposable graphs that are more complex than paths, with the help of constructions from~\cite{acharya2015string}.

\section{Number of labelings and subgraphs}
\label{SecLabels}

Reconstructing a non-confusable graph $G$ can be done by simply enumerating all its compositions and comparing those to the possible labelings of $G$ over an alphabet of size $k$. However, doing so requires the enumeration of those labelings. In this section, we will see that the number of such labelings is extremely high (even with the condition that they should be non-isomorphic), and as such this method is generally inefficient, justifying the specific analysis we do in later sections.
Furthermore, if we consider a labeling of a given graph as an encoding of information, then we are interested in graphs which have as many labelings as possible.  

\begin{definition}
	For a given graph $G$, we denote by $\rchi_k(G)$ the number of non-isomorphic labelings of $G$ over an alphabet of size $k$.
\end{definition}

We have included in this section proofs for the value $\rchi_k(G)$ over the graph classes studied in this article. Although some of these values may be known in the literature (such as $\rchi_k(P_n)$), we have decided to include the proofs as they are short and of independent interest, as well as for the sake of completeness.
We begin by giving the value of $\rchi_k(G)$ for paths and then continue with complete graphs.

\begin{theorem}\label{thmPathLabels}
	Given a path $P_n$ on $n$ vertices, we have $$\rchi_k(P_n)=\frac{k^n+k^{\lceil n/2\rceil}}{2}.$$
\end{theorem}
\begin{proof}
	Given a path $P_n$, the total number of labelings is $k^n$. However, each labeling is isomorphic with its reverse. Note that, while they are isomorphic, a reverse is different from the original labeling unless the original labeling is a palindrome. Furthermore, the $\lceil n/2\rceil$ first labels of a path fix the latter $\lfloor n/2\rfloor$ labels of a palindrome. Hence, we have $\rchi_k(P_n)=\frac{k^n-k^{\lceil n/2\rceil}}{2}+k^{\lceil n/2\rceil}=\frac{k^n+k^{\lceil n/2\rceil}}{2}$.
\end{proof}

\begin{theorem}\label{thmKnLabels}
	Given a complete graph $K_n$ on $n$ vertices, we have $\rchi_k(K_n)=\binom{n+k-1}{n}$.
\end{theorem}
\begin{proof}
	As no vertex is distinguishable from others, the problem can be translated as: in how many ways can we place $n$ balls into $k$ bins? Hence, by \cite{stanley2011enumerative}, the value $\rchi_k(K_n)$ can be solved using the ``stars and bars-technique'', which gives $\rchi_k(K_n)=\binom{n+k-1}{n}.$  
\end{proof}

We now determine the value of $\rchi_k(G)$ for the graph classes that we will study in \Cref{sec-reconstructing}.

\begin{theorem}\label{thmSubStarLabels}
	Given a subdivided star $S$ with maximum degree $\Delta \geq 3$, maximum branch length $m$, and $n_i$ branches of length $i$ for $i \in \{1,\ldots,m\}$, we have $\rchi_k(S)=k\prod_{i=1}^{m}\binom{k^i+n_i-1}{n_i}.$
\end{theorem}
\begin{proof}
	As we can identify the ends of a branch from each other (only one of them is a leaf), there are $k^i$ non-isomorphic labelings inside a single branch of length $i$. However, the $n_i$ branches of length $i$ are in symmetric positions and hence we can consider labeling them as labeling vertices of the complete graph $K_{n_i}$ using an alphabet containing $k^i$ symbols. Since we also have $k$ options for the center, by \Cref{thmKnLabels}, we have $$\rchi_k(S)=k\prod_{i=1}^{m}\rchi_{k^i}(K_{n_i})=k\prod_{i=1}^{m}\binom{k^i+n_i-1}{n_i}.\qedhere$$
\end{proof}

\begin{corollary}\label{corStarLabels}
	Given a star $S_n$ on $n+1$ vertices, we have $\rchi_k(S_n)=k\binom{n+k-1}{n}.$
\end{corollary}

\begin{corollary}\label{corStarS11mLabels}
	Given a subdivided star $S_{1,1,m}$ on $m+3$ vertices, we have $\rchi_k(S_{1,1,m})=\frac{k^{m+3}+k^{m+2}}{2}.$
\end{corollary}
\begin{proof}
	By \Cref{thmSubStarLabels}, we have $\rchi_k(S_{1,1,m})=k \binom{k+1}{2}\binom{k^m}{1}=\frac{k^{m+3}+k^{m+2}}{2}$.
\end{proof}

\begin{theorem}\label{thmBistarLabels}
	Given a bistar $B_{n,m}$, we have $$\rchi_k(B_{n,m})=\begin{cases}
		k^2\binom{n+k-1}{n}\binom{m+k-1}{m}, &\text{ if } n\neq m;\\
		\frac{k^2\binom{n+k-1}{n}^2+k\binom{n+k-1}{n}}{2}, &\text{ if } n= m.
	\end{cases}$$
\end{theorem}
\begin{proof}
	Let us first consider the case with $n\neq m$. Note that in this case the two stars around the  centers are not in symmetrical positions. Hence, we may choose both centers in $k$ ways, each leaf around one of the centers in $\binom{n+k-1}{n}$ ways and the leaves around the other center in $\binom{m+k-1}{m}$ ways. Hence, we have $\rchi_k(B_{n,m})=k^2\binom{n+k-1}{n}\binom{m+k-1}{m}$.
	
	Let us next consider the case with $n=m$. In this case, the two stars around the centers are in symmetrical positions. By \Cref{corStarLabels}, there are $k\binom{n+k-1}{n}$ ways to label one half of the bistar. As there are two symmetric halves in the bistar, we may consider labeling a path $P_2$ using an alphabet of size $k\binom{n+k-1}{n}$. Hence, $\rchi_k(B_{n,m})=\rchi_{k\binom{n+k-1}{n}}(P_2)=\frac{k^2\binom{n+k-1}{n}^2+k\binom{n+k-1}{n}}{2}$, by \Cref{thmPathLabels}.
\end{proof}

Note that for $n$-vertex graphs,  $\rchi_k(S_{1,1,n-3})=\frac{k^{n}+k^{n-1}}{2}$ by \Cref{corStarS11mLabels}  while \Cref{thmPathLabels} gives $\rchi_k(P_n)=\frac{k^n+k^{\lceil n/2\rceil}}{2}.$ This implies that we may store more information in $S_{1,1,n-3}$ than in $P_n$. When $k=2$ and $n$ is large, the ratio  $\rchi_k(S_{1,1,n-3}) : \rchi_k(P_n) \approx 1.5$.

\section{Reconstructability of some trees}
\label{sec-reconstructing}

In this section, we will show how to reconstruct several subclasses of trees. We begin with useful lemmas, then each subsection will be dedicated to a specific subclass.

\begin{lemma}
	\label{lem-leavesAndInternal}
	One can find the composition of leaves and internal vertices of a tree $\mathcal{T}$ using two sum-queries.
\end{lemma}

\begin{proof}
	Let $\mathcal{T}=(T,\lambda)$ be such that $T$ is a tree of order $n$ with $\ell$ leaves.
	The query $S_n(\mathcal{T})$ gives us the list of all labels. The query $S_{n-1}(\mathcal{T})$ gives us the sum of labels for all the subtrees of order $n-1$, each of which corresponds to $T$ minus exactly one leaf. Hence, $\ell$ subtrees are being summed in $S_{n-1}(\mathcal{T})$; each leaf is counted $\ell-1$ times and each internal vertex is counted $\ell$ times. Thus, the composition $I$ of internal vertices is given by $S_{n-1}(\mathcal{T})- (\ell-1) S_{n}(\mathcal{T})$, and the composition $L$ of leaves is given by $S_{n}(\mathcal{T}) - I$.
\end{proof}

\begin{example}
	Consider the labeled tree $\mathcal{T}$ on nine vertices presented in Figure \ref{fig:Ex1graph}. The queries $S_9(\mathcal{T})$ and $S_{8}(\mathcal{T})$ output $(3,1,2,3)$ and $(16,6,10,16)$ for alphabet $\sum=(A,B,C,D)$, respectively. As there are six leaves, the composition of the internal vertices is given by $I=(16,6,10,16)-(6-1)\times (3,1,2,3)=(1,1,0,1)$: the internal vertices have labels $A$, $B$ and $D$. From this, we can extract the composition of leaves: $S_9(\mathcal{T})-I=(3,1,2,3)-(1,1,0,1)=(2,0,2,2)$, so two leaves have label $A$, two leaves have label $C$, and two leaves have label $D$.
\end{example}

\begin{corollary}
	\label{cor-reconstructingStars}
	Stars are sum-reconstructable using two sum-queries.
\end{corollary}

\begin{lemma}
	\label{lem-centerOfAStar}
	The label of the center of a subdivided star with maximum degree $\Delta \geq 3$ can be found using three sum-queries.
\end{lemma}

\begin{proof}
	Let $\mathcal{S}=(S,\lambda)$ be such that $S$ is a subdivided star with maximum degree $\Delta\geq3$. Use \Cref{lem-leavesAndInternal} to obtain $L$ and $I$ the composition of leaves and internal vertices, respectively. By \Cref{lem-s2G}, the query $S_{2}(\mathcal{S})$ counts each leaf once, each degree~2 vertex twice, and the center $\Delta$ times. Hence, the label of the center (which is an internal vertex) is given by the only non-zero coordinate of $S_{2}(\mathcal{S}) - L - 2I$.
\end{proof}

\subsection{Stars subdivided at most once}

\begin{theorem}
	\label{thm-reconstructingStarsSubdividedAtMostOnce}
	Let $k=2$. Stars subdivided at most once are reconstructable with four multiset-queries.
\end{theorem}

\begin{proof}
	Let $\Sigma=(A,B)$ and $\mathcal{S}=(S,\lambda)$ be such that $S$ is a star subdivided at most once. Note that, if $S$ is a path, then it contains at most five vertices and we are done by \cite{acharya2015string}, hence assume that there is a vertex of degree at least~$3$.
	First, use \Cref{lem-leavesAndInternal}  to obtain the composition of leaves and internal vertices and the label of the center vertex with \Cref{lem-centerOfAStar}, for a total of three queries. Recall that we know, from $S_{n}(\mathcal{S})$, the multiset of all labels for vertices of $S$. The query $\mathcal{M}_{n-2}(\mathcal{S})$ gives the composition of all substars of order $n-2$. Now, $R = \bigsqcup_{X \in \M_{n-2}(\mathcal{S})} \left( S_{n}(\mathcal{S}) - X \right)$ gives the compositions of the following: all branches of length~2, and all pairs of leaves. Since we know the composition of leaves, we can remove all pairs of leaves from $R$ and obtain the branches of length~2.
	Since we know the composition of branches of length~2 and the label of the center, we can deduce the labels of leaves adjacent to the center.
	
	Now, all that remains is to decide, for each branch of length~2, which vertex is a leaf and which is internal. However, since we know the composition of internal vertices (from which we can remove the center), we can apply the following: for each branch with composition $A^2$ (resp. $B^2$), a leaf adjacent to an internal vertex with label $A$ (resp. $B$) has label $A$ (resp. $B$) too; the remaining branches have composition $AB$ and thus leaves adjacent to an internal vertex with label $A$ (resp. $B$) have label $B$ (resp. $A$).
	This concludes the reconstruction of $\mathcal{S}$.
\end{proof}

\begin{remark}
	We note that by \cite[Theorem 2.2.4]{bartha2021reconstructibility}, we can identify the structures of subdivided stars from each others based on the numbers of their connected subgraphs and the number of vertices in these subgraphs. This implies that the proof of \Cref{thm-reconstructingStarsSubdividedAtMostOnce} (without restriction to 4 multiset-queries) can also be done without knowing the exact structure of the subdivided star (although we need to know that the underlying graph we are considering is a subdivided star).
\end{remark}

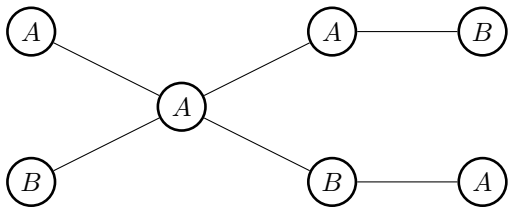
\begin{figure}[!h]
	\centering
	\begin{tikzpicture}
		\node[v] (c) at (2,1) {$A$};
		\node[v] (u1) at (4,2) {$A$};
		\node[v] (u2) at (6,2) {$B$};
		\node[v] (v1) at (4,0) {$B$};
		\node[v] (v2) at (6,0) {$A$};
		\node[v] (w1) at (0,2) {$A$};
		\node[v] (x1) at (0,0) {$B$};
		\foreach \I in {u,v,w,x} {\draw (c)to(\I1);}
		\draw (u1)to(u2);
		\draw (v1)to(v2);
	\end{tikzpicture}
	\caption{A labeled star subdivided at most once used for illustrating the proof of \Cref{thm-reconstructingStarsSubdividedAtMostOnce} on alphabet of size~2.}
	\label{fig-alphabetSize2IsNice}
\end{figure}

\medskip
\noindent\textbf{Example.} Consider the labeled star subdivided at most once $\mathcal{S}$ depicted in \Cref{fig-alphabetSize2IsNice}.

The queries $S_{7}(\mathcal{S})$ and $S_6(\mathcal{S})$ output $(4,3)$ and $(14,10)$, respectively. Since $(14,10)-(4-1) \times (4,3)=(2,1)$, two (resp. one) internal vertices have label $A$ (resp. $B$), and hence two leaves have label $A$ and two leaves have label $B$.
The query $S_{2}(\mathcal{S})$ outputs $(8,4)$. Since $S_{2}(\mathcal{S})-L-2I=(8,4) - (2,2) - 2 \times (2,1)=(2,0)$, the center vertex has label $A$, implying that the two degree~2 vertices have labels $A$ and $B$, respectively. Finally, the query $\mathcal{M}_{5}(\mathcal{S})$ outputs $\{\{ A^2B^3,6A^3B^2,A^4B^1\}\}$, hence $R = \bigsqcup_{X \in \mathcal{M}_{5}(\mathcal{S})} \left( S_{7}(\mathcal{S}) - X \right) = \{\{A^2, 6AB, B^2\}\}$. Since the combination of pairs of leaves is $\{\{A^2, 4AB, B^2\}\}$, the branches of length~2 both have a vertex labeled by $A$ and a vertex labeled by $B$. Removing $A^3B^2$ (for the center and those branches) from $A^4B^3$ gives the labels of the leaves adjacent to the center: $A$ and $B$.
As for the branches of length~2, they have composition $AB$, so the leaf adjacent to the internal vertex labeled $A$ (resp. $B$) has label $B$ (resp. $A$), and we are done.

\begin{remark}
	\label{rem-StarSubdividedAtMostOnce}
	Note that the reconstruction method described in the proof of \Cref{thm-reconstructingStarsSubdividedAtMostOnce} does not work for larger alphabets. This is due to the fact that the branches of length~2 offer more choices, as depicted in \Cref{fig-alphabetSize3IsBad}: we may know the composition of branches of length~2 and of leaves, but this does not allow for an immediate reconstruction. However, note that those two subdivided stars are not confusable, hence they are reconstructable but not sum-reconstructable.
	
	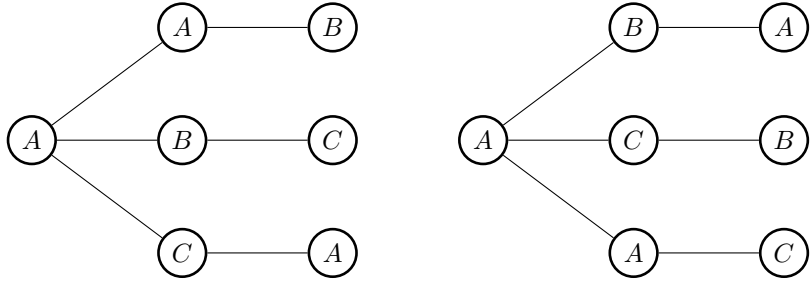
\begin{figure}[h]
		\centering
		\begin{tikzpicture}
			\node (conf1) at (0,0) {
				\begin{tikzpicture}
					\node[v] (c) at (0,1.5) {$A$};
					\node[v] (u1) at (2,0) {$C$};
					\node[v] (u2) at (4,0) {$A$};
					\node[v] (v1) at (2,1.5) {$B$};
					\node[v] (v2) at (4,1.5) {$C$};
					\node[v] (w1) at (2,3) {$A$};
					\node[v] (w2) at (4,3) {$B$};
					\foreach \I in {u,v,w} {
						\draw (c)to(\I1);
						\draw (\I1)to(\I2);
					}
				\end{tikzpicture}
			};
			
			\node (conf2) at (6,0) {
				\begin{tikzpicture}
					\node[v] (c) at (0,1.5) {$A$};
					\node[v] (u1) at (2,0) {$A$};
					\node[v] (u2) at (4,0) {$C$};
					\node[v] (v1) at (2,1.5) {$C$};
					\node[v] (v2) at (4,1.5) {$B$};
					\node[v] (w1) at (2,3) {$B$};
					\node[v] (w2) at (4,3) {$A$};
					\foreach \I in {u,v,w} {
						\draw (c)to(\I1);
						\draw (\I1)to(\I2);
					}
				\end{tikzpicture}
			};
		\end{tikzpicture}
		
		\caption{Two non-sum-reconstructable labeled stars subdivided at most once with alphabet size~3.}
		\label{fig-alphabetSize3IsBad}
	\end{figure}
\end{remark}

\subsection{Bistars}

\begin{theorem}
	\label{thm-reconstructingBistarsK2}
	Let $k=2$. Bistars are reconstructable with three multiset-queries.
\end{theorem}

\begin{proof}
	Let $\Sigma=(A,B)$ and $\mathcal{B}=(G,\lambda)$ be such that $G$ is a bistar with center vertices $u$ and $v$, vertex $u$ (resp. $v$) having $\ell_u$ (resp. $\ell_v$) adjacent leaves. First, use \Cref{lem-leavesAndInternal} to obtain the composition of leaves and of $u$ and $v$.
	Now, there are two cases to consider.
	
	\medskip\noindent\textbf{Case 1: $u$ and $v$ have the same label.} Assume without loss of generality that they both have label $A$. Let $a_u$ and $b_u$ (resp. $a_v$, $b_v$) be the (yet unknown) number of leaves attached to $u$ (resp. $v$) labeled with $A$ and $B$, respectively. Denote by $\alpha$ and $\beta$ the total number of vertices labeled with $A$ and $B$, respectively (which are known from $S_{n}(\mathcal{B})$). Note that we have $\alpha=a_u+a_v+2$ and $\beta=b_u+b_v$, as well as $a_u+b_u=\ell_u$ and $a_v+b_v=\ell_v$.
	
	The query $\mathcal{M}_{3}(\mathcal{B})$ outputs the composition of connected subgraphs of order~3, from which we can obtain $x$ and $y$ the number of subgraphs with compositions $A^3$ and $A^2B$, respectively. Note that if $x=0$ then we are done (all leaves have label $B$), and likewise if $y=0$ (all leaves have label $A$). Hence, assume that $x,y > 0$, we have:
	\begin{equation}
		\label{eqn-bistar}
		x = a_u + a_v + \binom{a_u}{2} + \binom{a_v}{2} = \frac{\alpha-2+a_u^2+a_v^2}{2} \Longrightarrow a_u^2+a_v^2 = 2x - \alpha +2
	\end{equation}
	as well as:
	\begin{align*}
		y & = b_u + b_v + a_ub_u + a_vb_v & \\
		& = \beta + a_u(\ell_u-a_u) + a_v(\ell_v - a_v) & \text{by $a_u+b_u=\ell_u$ and $a_v+b_v=\ell_v$} \\
		& = \beta + \ell_u a_u + \ell_v a_v - (a_u^2+a_v^2) & \\
		& = \beta + \ell_u a_u + \ell_v (\alpha - a_u - 2) - (a_u^2+a_v^2) & \text{by $\alpha=a_u+a_v+2$} \\
		& = \beta + \ell_u a_u + \ell_v (\alpha - a_u - 2) - (2x-\alpha+2) & \text{by (\ref{eqn-bistar})} \\
		& = \beta + a_u(\ell_u - \ell_v) + \alpha (\ell_v + 1) -2 (\ell_v+x+1) & 
	\end{align*}
	and hence:
	$$ a_u = \frac{y - \beta - \alpha (\ell_v+1) + 2(\ell_v+x+1)}{\ell_u-\ell_v} $$
	from which we can now find $a_v$, $b_u$ and $b_v$.

	\medskip\noindent\textbf{Case 2: $u$ and $v$ have different labels.} Let $a$ (resp. $y$) and $x$ (resp. $b$) be the (yet unknown) number of leaves attached to the internal vertex with label $A$ (resp. $B$) and themselves labeled with $A$ and $B$, respectively. Denote by $\alpha$ and $\beta$ the total number of vertices labeled with $A$ and $B$, respectively (which are known from $S_{n}(\mathcal{B})$). Note that we have $\alpha=a+y+1$ and $\beta=b+x+1$.
	
	The query $\mathcal{M}_{2}(\mathcal{B})$ outputs the composition of pairs of connected vertices, from which we can obtain $\alpha'$ and $\beta'$ the number of such pairs with composition $A^2$ and $B^2$, respectively. Note that we necessarily have $a = \alpha'$ and $b = \beta'$. This allows us to also obtain $x=\beta-\beta'-1$ and $y=\alpha-\alpha'-1$. We can now use those values to deduce whether $\ell_u=x+a$ and $\ell_v=b+y$, or $\ell_v=x+a$ and $\ell_u=y+b$; note that the two possibilities can  be true simultaneously if and only if $\ell_u = \ell_v$ (in which case, we are trivially done). Hence, we can find the labels of $u$ and $v$ as well as the labels of the attached leaves.
	
	\medskip
	This concludes the case analysis; note that we use two queries, followed by one more query in both cases, so three queries in total.
\end{proof}

\begin{theorem}
	\label{thm-reconstructingBistarsAnyK}
	For $k \geq 3$, bistars are reconstructable. Furthermore, if the two non-leaves have different labels, then, this can be done in three multiset-queries.
\end{theorem}

\begin{proof}
	Let $\mathcal{B}=(B,\lambda)$ be such that $B$ is a bistar with center vertices $u$ and $v$, vertex $u$ (resp. $v$) having $\ell_u$ (resp. $\ell_v$) adjacent leaves. First, use \Cref{lem-leavesAndInternal} to obtain the composition of leaves and of $u$ and $v$. There are two cases to consider.
	
	\medskip\noindent\textbf{Case 1: $u$ and $v$ have different labels.} We will apply a similar technique as in the proof of the same case in \Cref{thm-reconstructingBistarsK2}. Assume without loss of generality that $u$ and $v$ have labels $X$ and $Y$, denote the (yet unknown) number of leaves attached to the vertex with label $X$ (resp. $Y$) by $\ell_x$ (resp. $\ell_y$); we want to decide whether $\ell_u=\ell_x$ or $\ell_u=\ell_y$. Note that we also know $x$ (resp. $y$) the number of vertices labeled with $X$ (resp. $Y$).
	
	The query $\mathcal{M}_2(\mathcal{B})$ outputs the composition of pairs of connected vertices. There are three types of combinations: $XY$, $XZ$ for $Z \neq Y$ and $YZ$ for $Z \neq X$. The combinations $XZ$ for $Z \neq Y$ give us all leaves adjacent to the vertex with label $X$, denote by $z_x$ their number. The combinations $YZ$ for $Z \neq X$ give us all leaves adjacent to the vertex with label $Y$, denote by $z_y$ their number. From those, we also know $x_2$ and $y_2$, the numbers of combinations $XX$ and $YY$, respectively. Finally, let $\alpha$ (resp. $\beta$) be the (yet unknown) number of vertices labeled with $Y$ (resp. $X$) attached to the internal vertex with label $X$ (resp. $Y$). We have $x=\beta+x_2+1$ and $y=\alpha+y_2+1$, from which we can deduce $\alpha$ and $\beta$. Furthermore, $\ell_x=z_x+\alpha$ and $\ell_y=z_y+\beta$, so we can decide whether $\ell_u=\ell_x$ or $\ell_u=\ell_y$ (if both are true simultaneously, then, $\ell_u=\ell_v$ and we are trivially done).
	
	\medskip\noindent\textbf{Case 2: $u$ and $v$ have the same label.}
	
	Let us denote by $A$ the label of vertices $u$ and $v$. Since $u$ and $v$ have the same label, we may consider the largest composition which contains only one $A$ (this is found by exploring queries $\mathcal{M}_i(\mathcal{B})$ for decreasing values of $i$, starting from $i=\ell_u+\ell_v+2$). The single $A$ in this composition is either $u$ or $v$ (unless every vertex is labeled with $A$, in which case the reconstruction is trivially done). Hence, this composition gives every non-$A$ label on one side of the bistar. Denote the total number of these non-$A$ symbols by $\beta_1$. As we can obtain the total number of each symbol with 1-vertex subgraphs, we can also obtain the number of non-$A$ symbols on the other side of the bistar, denote this by $\beta_2$. Hence, we only need to solve whether we have $\beta_1$ non-$A$ symbols adjacent to $u$ or adjacent to $v$. Note that when $\beta_1=\beta_2$ or $\ell_u=\ell_v$, we are already done. 
	
	Assume without loss of generality that $\ell_v> \ell_u$ and let $\ell_v=\ell_u+\ell$. We denote by $A_u$ ($A_v$) and $B_u$ ($B_v$) the number of symbols $A$ and non-$A$ symbols adjacent to $u$ ($v$), respectively.
	Consider next the number of compositions $A^3$ (obtained using $\mathcal{M}_3(\mathcal{B})$) and denote it by $a_3$. Since each $A^3$ contains either two vertices adjacent to one of $u$ or $v$, or both $u$ and $v$ and a single vertex adjacent to one of $u$ or $v$, we have $$a_3=\binom{A_u}{2}+\binom{A_v}{2}+A_u+A_v=\frac{A_u^2+A_v^2+A_u+A_v}{2}.$$ Furthermore, as we know the total number of vertices labeled with $A$, denote it by $\alpha$, we also know that $A_u+A_v=\alpha-2$. Hence, $2a_3-\alpha+2=A_u^2+A_v^2$. Furthermore, we have $$A_u^2+A_v^2=(\ell_u-B_u)^2+(\ell_v-B_v)^2=\ell_u^2-2\ell_uB_u+B_u^2+\ell_v^2-2\ell_vB_v+B_v^2.$$ Note that we have $\beta_1^2+\beta_2^2=B_u^2+B_v^2$. Hence, we have $\frac{-2a_3+\alpha-2+\beta_1^2+\beta_2^2+\ell_u^2+\ell_v^2}{2}=\ell_uB_u+\ell_vB_v=(B_u+B_v)\ell_u+\ell B_v$. Therefore, $$B_v=\frac{-2a_3+\alpha-2+\beta_1^2+\beta_2^2+\ell_u^2+\ell_v^2-2(\beta_1+\beta_2)\ell_u}{2\ell}.$$
	As we can compute everything on the right hand side, we can compute $B_v$. Furthermore, $\{B_v,B_u\}=\{\beta_1,\beta_2\}$. As we mentioned above, this is enough to reconstruct the labeling of $\mathcal{B}$.
\end{proof}

\subsection{Stars $S_{1,1, m}$}

In this section, we study paths with one leaf false-twinned (that is, we create a false twin of one leaf). As shown with the proof of \Cref{thm-reconstructingS11k}, this very simple operation changes the status of the graph: recall that some paths are confusable~\cite{acharya2015string}; we show that false-twinning one of the leaves makes the graph sum-reconstructable with a linear number of queries, and the proof is quite technical. However, we will see later (\Cref{thm-triangleWithTail} in \Cref{secSmallConfusable}) that true-twinning a leaf (that is, creating a twin of it) creates a graph class exhibiting another behaviour (alternating between sum-reconstructable and confusable with the parity of $m$).

Recall that we have shown in \Cref{SecLabels} that subdivided stars of type $S_{1,1,m}$ allow more labelings than paths using $m+3$ labels. Furthermore, since it is relatively simple to construct $S_{1,1,m}$, these graphs seem suitable candidates for applications related to advanced information storages.

\begin{theorem}
	\label{thm-reconstructingS11k}
	For any nonnegative integer $m$, the subdivided star $\mathcal{S}=(S_{1,1, m},\lambda)$ is sum-reconstructable.
\end{theorem}

\begin{proof}
	Note that $P_3$ is sum-reconstructable by \Cref{cor-reconstructingStars}, so let $m$ be a positive integer.
	Let $\mathcal{S}=(S,\lambda)$ where $S=S_{1,1, m}$, denote the two leaves attached to the center of $S$ by $u_1$ and $u_2$, and the vertices from the center to the end of the branch of length $m$ by $v_1,v_2,\ldots,v_{m+1}$ (so the center is $v_1$ and $v_{m+1}$ is a leaf).
	
	For readability, we will use the notation of a given query $S_t(\mathcal{S})$ as a sum of labels, such that $X$ corresponds to the label of vertex $x$. For example, \Cref{lem-s2G} implies that for every positive $m$:
	\begin{equation}
		\label{eqn-S2S}
		S_2(\mathcal{S})=U_1+U_2+V_{m+1}+3V_1+2(V_2+\ldots+V_m).
	\end{equation}
	
	Note that, from \Cref{lem-centerOfAStar}, we can compute the label of $v_1$, and from \Cref{lem-leavesAndInternal} the composition of the three leaves ($U_1+U_2+V_{m+1}$) and of all internal vertices ($\sum_{i=2}^{m}V_i$). Furthermore, \Cref{cor-reconstructingStars} implies that the case $m \leq 1$ already holds. First, let $m = 2$. Since we already know the label of $v_1$, the other internal label is $v_2$. Furthermore, we have $S_3(\mathcal{S}) = 2U_1+2U_2+4V_1+3V_2+V_3$, allowing us to reconstruct $u_1$ and $u_2$, and thus $v_3$, completing the reconstruction of $\mathcal{S}$. Hence, for the rest of the proof, let us assume that $m \geq 3$.
	
	There are three types of subgraphs of order $i$ in $\mathcal{S}$ (for $3 \leq i \leq m+1$): one subgraph containing both $u_1$ and $u_2$ as well as $v_1,\ldots,v_{i-2}$, two subgraphs each containing one of either $u_1$ or $u_2$ as well as $v_1,\ldots,v_{i-1}$, and $m+2-i$ paths containing vertices $v_j,\ldots,v_{j+i-1}$ for $j \in \{1,\ldots,m+2-i\}$.
	
	In the following, we consider two cases, depending on the parity of $m$.
	
	\medskip\noindent\textbf{Case 1: $m$ is odd, so $m=2p-1$ for $p \geq 3$.} We will consider $2p-2$ queries, which will give the following results.
	First, for $i \in \{3,\ldots,p\}$,
	\begin{equation}\label{eqSi}
		S_i(\mathcal{S}) = 2(U_1+U_2) + \sum_{j=1}^{i-2} \left((j+3)V_j\right) + (i+1)V_{i-1} + \sum_{j=i}^{2p-i+1} \left( iV_j \right) + \sum_{j=1}^{i-1} \left( jV_{2p-j+1} \right). 
	\end{equation}
	Then, for $i \in \{1,\ldots,p\}$,
	\begin{equation}\label{eqSpi}
		S_{p+i}(\mathcal{S}) = 2(U_1+U_2) + \sum_{j=1}^{p-i} \left( (j+3)V_j \right) + \sum_{j=p-i+1}^{p+i-2} \left( (p-i+4)V_j \right) + (p-i+3)V_{p+i-1} + \sum_{j=1}^{p-i+1} \left( jV_{2p-j+1} \right).     
	\end{equation}

	\begin{claim}
		\label{clm-oddBranch-1}
		For $i \in \{1,\ldots,p-1\}$, \[S_{p+i}(\mathcal{S}) - S_{p+i+1}(\mathcal{S}) = \sum_{j=p-i+1}^{p+i-2} (V_j) - V_{p+i}.\]
	\end{claim}
	
	\begin{proof}[Proof of Claim \ref{clm-oddBranch-1}]
		By applying Equation (\ref{eqSpi}) twice, we obtain:
		\begin{eqnarray*}
			S_{p+i}(\mathcal{S}) - S_{p+i+1}(\mathcal{S}) & = & \left(2(U_1+U_2) + \sum_{j=1}^{p-i} \left( (j+3)V_j \right) + \sum_{j=p-i+1}^{p+i-2} \left( (p-i+4)V_j \right)\right. \\
			& & \left.+ (p-i+3)V_{p+i-1} + \sum_{j=1}^{p-i+1} \left( jV_{2p-j+1} \right)\right) \\
			& & - \left(2(U_1+U_2) + \sum_{j=1}^{p-i-1} \left( (j+3)V_j \right) + \sum_{j=p-i}^{p+i-1} \left( (p-i+3)V_j \right)\right. \\
			& & \left.+ (p-i+2)V_{p+i} + \sum_{j=1}^{p-i} \left( jV_{2p-j+1} \right)\right) \\
			& = & (p-i+3 - (p-i+3)) V_{p-i} + \sum_{j=p-i+1}^{p+i-2} \left( V_j \right) \\
			& & + (p-i+3 - (p-i+3))V_{p+i-1} + (p-i+1-(p-i+2))V_{p+i} \\
			& = & \sum_{j=p-i+1}^{p+i-2} (V_j) - V_{p+i}.\qedhere
		\end{eqnarray*}
	\end{proof}
	
	\begin{claim}
		\label{clm-oddBranch-2}
		For $i \in \{3,\ldots,p\}$, \[S_{i+1}(\mathcal{S}) - S_{i}(\mathcal{S}) = V_{i-1} + 2 V_i + \sum_{j=i+1}^{2p-i} (V_j).\]
	\end{claim}
	
	\begin{proof}[Proof of Claim \ref{clm-oddBranch-2}]
		By applying Equation (\ref{eqSi}) twice for $i\in\{3,\ldots, p-1\}$, we obtain:
		\begin{eqnarray*}
			S_{i+1}(\mathcal{S}) - S_{i}(\mathcal{S}) & = & \left(2(U_1+U_2) + \sum_{j=1}^{i-1} \left( (j+3)V_j \right) + (i+2)V_i\right. \\
			& &\left. + \sum_{j=i+1}^{2p-i} \left( (i+1)V_j \right) + \sum_{j=1}^{i} \left( jV_{2p-j+1} \right)\right) \\
			& & - \left(2(U_1+U_2) + \sum_{j=1}^{i-2} \left( (j+3)V_j \right) + (i+1)V_{i-1} \right.\\
			& &\left. + \sum_{j=i}^{2p-i+1} \left( iV_j \right) + \sum_{j=1}^{i-1} \left( jV_{2p-j+1} \right)\right) \\
			& = & (i+2 - (i+1)) V_{i-1} + (i+2 - i)V_i + \sum_{j=i+1}^{2p-i} \left( V_j \right) + (i-i)V_{2p-i+1} \\
			& = & V_{i-1} + 2 V_i + \sum_{j=i+1}^{2p-i} (V_j).
		\end{eqnarray*}
		
		We note that $S_{p+1}(\mathcal{S})-S_p(\mathcal{S})$ can be obtained with similar calculations by applying (\ref{eqSpi}) for $S_{p+1}(\mathcal{S})$.
	\end{proof}
	
	We next compute the labels of vertices in the following way: alternating between $S_{p+i+1}(\mathcal{S})-S_{p+i}(\mathcal{S})$ for $i \in \{1,\ldots,p-1\}$ and $S_{i'+1}(\mathcal{S})-S_{i'}(\mathcal{S})$ for $i' \in \{p,\ldots,3\}$, we will be able to find the labels of $v_{p+i}$ and both $v_{i'-1}$ and $v_{i'}$, respectively. We prove that we can do those computations by induction on $i \in \{1,\ldots,p-2\}$ and $i' \in \{p,\ldots,3\}$ (with $i$ increasing and $i'$ decreasing at each step; so $i+i'=p+1$). First, if $(i,i')=(1,p)$ and $p>2$, we compute, by \Cref{clm-oddBranch-1}:
	\[ S_{p+2}(\mathcal{S}) - S_{p+1}(\mathcal{S}) = V_{p+1}-\sum_{j=p-1+1}^{p+1-2} (V_j) = V_{p+1} \]
	which gives the label of $v_{p+1}$, and, by \Cref{clm-oddBranch-2}:
	\[ S_{p+1}(\mathcal{S}) - S_{p}(\mathcal{S}) = V_{p-1} + 2 V_p + \sum_{j=p+1}^{2p-p} (V_j) = V_{p-1} + 2V_p \]
	which gives the labels of $v_{p-1}$ and $v_p$ (if $V_p=V_{p-1}$, the same label appears thrice; otherwise, $v_p$ has the label that appears twice and $v_{p-1}$ the label that appears once).
	Note that, if $p=2$ (so we are reconstructing $\mathcal{S}=(S_{1,1,3},\lambda)$), there is an adjustment: we apply \Cref{clm-oddBranch-1} with $i=1$ and compute $S_4(\mathcal{S})-S_3(\mathcal{S})=V_3$. We can then deduce the label of $v_2$ which is the last of the three internal vertices (recall that we have already obtained the label of the center vertex $v_1$). Next, we use Equations~(\ref{eqn-S2S}) and (\ref{eqSi}) to compute $S_3(\mathcal{S})-S_2(\mathcal{S})=U_1+U_2+V_1+2V_2$, obtaining the labels of $u_1$ and $u_2$, after which we can find the label of $v_4$. This ends the reconstruction for $p=2$. For the rest of Case~1, we assume that $p \geq 3$.
	
	Assume now that the statement holds for a pair $(i-1,i'+1)$, so we know the labels of $v_{p-i+1},\ldots,v_{p+i-1}$. We compute, by \Cref{clm-oddBranch-1}:
	\begin{eqnarray*}
		S_{p+i+1}(\mathcal{S}) - S_{p+i}(\mathcal{S}) & = & V_{p+i} - \sum_{j=p-i+1}^{p+i-2} (V_j) \\
		& = & V_{p+i} - \sum_{j=i'}^{p+i-2} (V_j) \\
		& = & V_{p+i} - Q
	\end{eqnarray*}
	where $Q$ is a known quantity by induction hypothesis, hence we obtain the label of $v_{p+i}$. Therefore, we obtain labels $V_{p},\ldots, V_{p+i}$.
	Using now \Cref{clm-oddBranch-2}, we compute:
	\begin{eqnarray*}
		S_{i'+1}(\mathcal{S}) - S_{i'}(\mathcal{S}) & = & V_{i'-1} + 2 V_{i'} + \sum_{j=i'+1}^{2p-i'} (V_j) \\
		& = & V_{i'-1} + 2 V_{i'} + \sum_{j=i'+1}^{p+i-1} (V_j) \\
		& = & V_{i'-1} + Q' \\
	\end{eqnarray*}
	where $Q'$ is a known quantity by induction hypothesis, hence we obtain the label of $v_{i'-1}$. Hence, when $i$ traverses from $1$ to $p-2$ and $i'$ from $p$ to $3$, we are able to compute the labels of $v_2,\ldots,v_{2p-2}$ (and the label of $v_1$ was already known). Now, using one last time \Cref{clm-oddBranch-1}, we compute
	\[ S_{2p}(\mathcal{S}) - S_{2p-1}(\mathcal{S}) = V_{2p-1} - \sum_{j=2}^{2p-3} (V_j) = V_{2p-1} + Q'' \]
	where $Q''$ is a known quantity, allowing us to find the label of $v_{2p-1}$.
	We now only have to find the label of $v_{2p}$. We observe that $S_{2p+1}(\mathcal{S}) =  2(U_1+U_2) + \sum_{j=1}^{2p-1} \left( 3V_j \right) + 2V_{2p}$, which leads to:
	\begin{eqnarray*}
		S_{2p+1}(\mathcal{S}) - S_{2p}(\mathcal{S}) & = & 2(U_1+U_2) + \sum_{j=1}^{2p-1} \left( 3V_j \right) + 2V_{2p} - 2(U_1+U_2) - \sum_{j=1}^{2p-2} \left( 4V_j \right) - 3V_{2p-1} - V_{2p} \\
		& = & V_{2p} - \sum_{j=1}^{2p-2} V_j \\
		& = & V_{2p} + Q'''
	\end{eqnarray*}
	where $Q'''$ is a known quantity. The labels of $u_1$ and $u_2$ are now the two labels of leaves remaining. Recall that we already found the composition $U_1+U_2+V_{m+1}$. Hence, we can also deduce the labels of $u_1$ and $u_2$.
	This concludes the case where $m$ is odd.
	
	\medskip\noindent\textbf{Case 2: $m$ is even, so $m=2p$ for $p \geq 2$.}
	Let us first consider the case with $p=2$ (so we are reconstructing $\mathcal{S}=(S_{1,1,4},\lambda)$), which will serve as a concrete example of the case disjunction that we are going to need for this case. This is necessary since one main equation needed for the proof (\Cref{clm-evenBranch-1}) does not hold for this small graph. We have $S_5(\mathcal{S})-S_4(\mathcal{S})=V_4-V_2$. We divide our considerations according to whether we have $V_2 \neq V_4$ or $V_2 = V_4$. If $V_2\neq V_4$, then we can know $V_2$ and $V_4$. Hence, we may deduce $V_3=I-V_1-V_2-V_4$ and $U_1,U_2$ from $S_3(\mathcal{S})-S_2(\mathcal{S})=U_1+U_2+V_1+2V_2+V_3$. Finally, we may finish with $V_5=L-U_1-U_2$.
	If $V_2=V_4$, then we compute $S_4(\mathcal{S})-S_3(\mathcal{S})=V_2+V_3$, so if $V_2=V_3$ we are done (this gives their value, and we then conclude as above). Assume then that $V_2 \neq V_3$, we use Equation~(\ref{eqn-S2S}) to compute $S_6(\mathcal{S})-S_2(\mathcal{S})=U_1+U_2+V_2+V_3+V_4+V_5=L+2V_2+V_3$ where $L$ is the sum of the labels of leaves (which is known). This gives us the values of $V_2$ and $V_3$, and we then conclude as above.
	For the rest of Case~2, we assume that $p \geq 3$.
	
	We will consider $2p-1$ queries, which will give the following results. First, for $i \in \{3,\ldots,p+1\}$,
	\begin{equation}\label{eqSiEven}
		S_i(\mathcal{S}) = 2(U_1+U_2) + \sum_{j=1}^{i-2} \left( (j+3)V_j \right) + (i+1)V_{i-1} + \sum_{j=i}^{2p-i+2} \left( iV_j \right) + \sum_{j=1}^{i-1} \left( jV_{2p-j+2} \right). 
	\end{equation}
	Then, for $i \in \{1,\ldots,p\}$,
	\begin{equation}\label{eqSpiEven}
		S_{p+1+i}(\mathcal{S}) = 2(U_1+U_2) + \sum_{j=1}^{p-i+1} \left( (j+3)V_j \right) + \sum_{j=p-i+2}^{p+i-1} \left( (p+4-i)V_j \right) + (p+3-i)V_{p+i} + \sum_{j=1}^{p-i+1} \left( jV_{2p-j+2} \right). \end{equation}
	
	\begin{claim}
		\label{clm-evenBranch-1}
		For $i \in \{3,\ldots,p\}$:
		\[ S_{i+1}(\mathcal{S}) - S_i(\mathcal{S}) = V_{i-1} + 2V_i + \sum_{j=i+1}^{2p-i+1} \left( V_j \right). \]
	\end{claim}
	
	\begin{proof}[Proof of Claim \ref{clm-evenBranch-1}]        By applying Equation (\ref{eqSiEven}) twice, we obtain for $i \in \{3,\ldots,p\}$
		\begin{eqnarray*}
			S_{i+1}(\mathcal{S}) - S_i(\mathcal{S}) & = & \left(2(U_1+U_2) + \sum_{j=1}^{i-1} \left( (j+3)V_j \right) + (i+2)V_{i} + \sum_{j=i+1}^{2p-i+1} \left( (i+1)V_j \right) + \sum_{j=1}^{i} \left( jV_{2p-j+2} \right)\right) \\
			& & - \left(2(U_1+U_2) + \sum_{j=1}^{i-2} \left( (j+3)V_j \right) + (i+1)V_{i-1} + \sum_{j=i}^{2p-i+2} \left( iV_j \right) + \sum_{j=1}^{i-1} \left( jV_{2p-j+2} \right)\right) \\
			& = & \left( (i-1+3)-(i+1) \right) V_{i-1} + \left( i+2-i \right) V_i + \sum_{j=i+1}^{2p-i+1} \left( V_j \right) + (i-i)V_{2p-i+2} \\
			& = & V_{i-1} + 2V_i + \sum_{j=i+1}^{2p-i+1} \left( V_j \right).\qedhere
		\end{eqnarray*}
	\end{proof}
	
	\begin{claim}
		\label{clm-evenBranch-2}
		For $i \in \{ 1,\ldots,p-1\}$:
		\[ S_{p+1+i}(\mathcal{S}) - S_{p+2+i}(\mathcal{S}) = V_{p-i+1} - V_{p+i+1} + \sum_{j=p-i+2}^{p+i-1} \left( V_j \right). \]
	\end{claim}
	
	\begin{proof}[Proof of Claim \ref{clm-evenBranch-2}]  By applying Equation (\ref{eqSpiEven}) twice, we obtain for $i \in \{ 1,\ldots,p-1\}$
		\begin{eqnarray*}
			S_{p+1+i}(\mathcal{S}) - S_{p+2+i}(\mathcal{S}) & = & \left(2(U_1+U_2) + \sum_{j=1}^{p-i+1} \left( (j+3)V_j \right) + \sum_{j=p-i+2}^{p+i-1} \left( (p+4-i)V_j \right)\right. \\
			& & \left.+ (p+3-i)V_{p+i} + \sum_{j=1}^{p-i+1} \left( jV_{2p-j+2} \right)\right) \\
			& & - \left(2(U_1+U_2) + \sum_{j=1}^{p-i} \left( (j+3)V_j \right) + \sum_{j=p-i+1}^{p+i} \left( (p+3-i)V_j \right)\right. \\
			& &\left. + (p+2-i)V_{p+i+1} + \sum_{j=1}^{p-i} \left( jV_{2p-j+2} \right)\right) \\
			& = & \left( (p-i+1+3) - (p+3-i) \right) V_{p-i+1} + \sum_{j=p-i+2}^{p+i-1} \left( V_j \right) \\
			& & + \left( (p+3-i) - (p+3-i) \right) V_{p+i} + \left( (p-i+1) - (p+2-i) \right) V_{p+i+1} \\
			& = & V_{p-i+1} - V_{p+i+1} + \sum_{j=p-i+2}^{p+i-1} \left( V_j \right).\qedhere
		\end{eqnarray*}
	\end{proof}
	
	We now start computing the labels of vertices. First, note that by Claim \ref{clm-evenBranch-1} with $i=p$ we have:
	\begin{equation}\label{eqSp1Sp}
		S_{p+1}(\mathcal{S}) - S_p(\mathcal{S}) = V_{p-1} + 2V_p + V_{p+1}.
	\end{equation}
	
	We can derive two subcases from Equation~(\ref{eqSp1Sp}).
	
	\medskip\noindent\textbf{Subcase 2.1: $V_{p-1} \neq V_{p+1}$, or $V_{p-1}=V_p=V_{p+1}$.} In this case, $S_{p+1}(\mathcal{S}) - S_p(\mathcal{S})$ clearly allows us to deduce $V_p$ (it is the value that is repeated at least twice). From Equations (\ref{eqSiEven}) and (\ref{eqSpiEven}) for $i=p+1$ and $i=1$, respectively, we obtain $S_{p+2}(\mathcal{S})-S_{p+1}(\mathcal{S}) = V_p + V_{p+1}$. Since we already know $V_p$, we obtain the label of $v_{p+1}$, and thus the label of $v_{p-1}$ from Equation (\ref{eqSp1Sp}). Now, like in Case~1, we use induction and alternate between $S_{p+1+i}(\mathcal{S}) - S_{p+2+i}(\mathcal{S})$ and $S_{i'+1}(\mathcal{S})-S_{i'}(\mathcal{S})$, for increasing $i \in \{1,\ldots,p-3\}$ and decreasing $i' \in \{p-1,\ldots,3\}$, respectively.
	Doing so allows us to know the labels of $v_{i'-1},\ldots,v_{p+i+1}$ 
	for the pair $(i,i')$ where $i+i'=p$. Hence, our induction hypothesis is that after computing the pair $(i,i')$, we know the labels of $v_{i'-1},\ldots,v_{p+i+1}$. Once this is done, we finish with $i\in\{p-2,p-1\}$. By \Cref{clm-evenBranch-2} and induction hypothesis, we have
	\begin{equation}\label{eqC2.1eqi}
		S_{p+1+i}(\mathcal{S}) - S_{p+2+i}(\mathcal{S})=V_{p-i+1}-V_{p+i+1}+\sum_{j=p-i+2}^{p+i-1}(V_j),
	\end{equation}
	which allows us to compute the label of $v_{p+i+1}$. By \Cref{clm-evenBranch-1}, we also have
	$$S_{i'+1}(\mathcal{S})-S_{i'}(\mathcal{S})=V_{i'-1}+2V_{i'}+\sum_{j=i'+1}^{2p-i'+1}(V_j).$$
	Hence, the induction hypothesis together with the label of $v_{p+i+1}$ allows us to compute the label of $v_{i'-1}$. Note that as we have mentioned above, there are more potential values for $i$ than for $i'$ in order to compute all the internal labels. Hence, we first alternate between $i$ and $i'$ until the case $(p-2,3)$. At this point, we know the labels of $v_1,\dots, v_{2p-2}$. After this, we use Equation~(\ref{eqC2.1eqi}) for $i=p-2$ and $i=p-1$ which gives the labels of $v_{2p-1}$ and $v_{2p}$.
	We can again conclude by computing the label of $v_{2p+1}$ (from $S_{2p+1}(\mathcal{S})-S_{2p+2}(\mathcal{S})=\sum_{i=1}^{2p-1} (V_i) - V_{2p+1}$), allowing us to deduce the labels of $u_1$ and $u_2$, and thus completing the reconstruction.
	
	\medskip\noindent\textbf{Subcase 2.2: $V_{p-1} = V_{p+1}$ and $V_{p-1} \neq V_p$.} In this case, \Cref{eqSp1Sp} gives $S_{p+1}(\mathcal{S}) - S_p(\mathcal{S})= V_{p-1} + 2V_p + V_{p+1}$ and does not allow us to decide on the values of $V_{p-1}$ and $V_p$. Now, like in Case~1 and Subcase~2.1, we will use induction and alternate operations. However, this time, we will prove by induction either that, for values of $i$ that will be given later, $V_{2i}=V_{2i+2}$ and $V_{2i-1}=V_{2i+1}$ (depending on the parity of the induction step) and that $V_{2i-1}\neq V_{2i}$, or that we are able to determine the values of $V_{p-i},\ldots,V_{p+i}$ (after which we start a second induction, exactly like in Case~1 and Subcase~2.1, which will allow us to compute the remaining labels $V_{j}$).
	
	To start the induction, note that $V_{p-1}=V_{p+1}\neq V_p$, and let us denote for readability and without loss of generality that $V_{p-1}=V_{p+1}=A$ and $V_p=B$ (though we do not know these values yet). Assume that, for a given $0\leq i\leq \frac{p-3}{2}$, we have $V_{p-2j}=V_{p+2j}=B$ and $V_{p-2j-1}=V_{p+2j+1}=A$ for each $j \in \{0,\ldots,i\}$. This is our induction hypothesis, which holds for the base case $i=0$.
	
	Every main step of the induction is comprised of an \emph{even} and an \emph{odd} step, each seeing two computations using \Cref{clm-evenBranch-1,clm-evenBranch-2}. At any point during those two steps, we may actually collapse all values (that is, be able to know the values of $V_{p-1}$ and $V_p$, and thus all the other values that are equal to those) and end this induction. During the even step, we try to compute the values of $V_{p+2i}$ and then of $V_{p-2i}$. If we are able to determine their actual value, then the induction ends. Otherwise, we get $V_{p-2i}=V_{p+2i}=V_p$. Then, during the odd step, we try to compute the values of $V_{p+2i+1}$ and then of $V_{p-2i-1}$. Again, either we are able to determine them and the induction ends, or we get $V_{p-2i-1}=V_{p+2i+1}=V_{p-1}$. In the latter case, we continue the induction.
	
	\smallskip
	We first detail the \emph{even} case of the induction step, where we try to compute the values of $V_{p+2i}$ and then of $V_{p-2i}$. By induction hypothesis, we have $V_{p-2i+2}=V_{p-2i+4}=\cdots=V_{p+2i-2}=B$ and $V_{p-2i+1}=V_{p-2i+3}=\cdots=V_{p+2i-1}=A$.
	We compute, by \Cref{clm-evenBranch-2} and induction hypothesis (we have $1 \leq i \leq  \frac{p-3}{2}$):
	\begin{eqnarray*}
		S_{p+2i}(\mathcal{S}) - S_{p+1+2i}(\mathcal{S}) & = & V_{p-2i+2} - V_{p+2i} + \sum_{j=p-2i+3}^{p+2i-2} \left( V_j \right) \\
		& = & V_{p-2i+2} - V_{p+2i} + \sum_{\substack{j=p-2i+4 \\ j-p \text{ even}}}^{p+2i-2} \left( V_j \right) + \sum_{\substack{j=p-2i+3 \\ j-p \text{ odd}}}^{p+2i-3} \left( V_j \right) \\
		& = & B - V_{p+2i} + (2i-2)B + (2i-2)A.
	\end{eqnarray*}
	Now, if $V_{p+2i}=B$, then, we are still in the induction loop, since we cannot decide the value of $B$ with the equation. Indeed, we can in this case only know that $V_{p+2i}=V_{p}$. Otherwise, the equation allows us to ``collapse'' the unknown values of $V_{p-2i+1},\ldots,V_{p+2i-1}$, as well as the value of $V_{p+2i}$, in which case we end the induction. Recall that by induction hypothesis, we have $V_{p-2i+2}=V_{p-2i+4}=\cdots=V_{p+2i-2}=B$ and $V_{p-2i+1}=V_{p-2i+3}=\cdots=V_{p+2i-1}=A$.
	Similarly, using \Cref{clm-evenBranch-1}, the induction hypothesis, and the fact that $V_{p+2i}=B$ (we have $1\leq i \leq  \frac{p-3}{2}$):
	\begin{eqnarray*}
		S_{p+2-2i}(\mathcal{S}) - S_{p+1-2i}(\mathcal{S}) & = & V_{p-2i} + 2V_{p-2i+1} + \sum_{j=p-2i+2}^{p+2i} \left( V_j \right) \\
		& = & V_{p-2i} + 2A + (2i-1)A + (2i)B.
	\end{eqnarray*}
	Again, if $V_{p-2i} = B$, then, we cannot decide the values with the equation and thus, we remain in the induction loop; otherwise, we can deduce the value of $B$ as well as the unknown values of $V_{2p-2i+1},\ldots,V_{p+2i-1}$. In the latter case, we end this induction.
	
	\smallskip
	Let us next detail the \emph{odd} case of the induction step, where we try to compute the values of $V_{p+2i+1}$, and then of $V_{p-2i-1}$. Note that it is very similar to the even case.
	By induction hypothesis and after the even step, we have $V_{p-2i}=V_{p-2i+2}=\cdots=V_{p+2i}=B$ and $V_{p-2i+1}=V_{p-2i+3}=\cdots=V_{p+2i-1}=A$.
	We compute, by \Cref{clm-evenBranch-2} and induction hypothesis (we have $1 \leq i \leq  \frac{p-3}{2}$):
	\begin{eqnarray*}
		S_{p+1+2i}(\mathcal{S}) - S_{p+2+2i}(\mathcal{S}) & = & V_{p-2i+1} - V_{p+2i+1} + \sum_{j=p-2i+2}^{p+2i-1} \left( V_j \right) \\
		& = & V_{p-2i+1} - V_{p+2i+1} + \sum_{\substack{j=p-2i+2 \\ j-p \text{ even}}}^{p+2i-2} \left( V_j \right) + \sum_{\substack{j=p-2i+3 \\ j-p \text{ odd}}}^{p+2i-1} \left( V_j \right) \\
		& = & A - V_{p+2i+1} + (2i-1)B + (2i-1)A.
	\end{eqnarray*}
	Again, if $V_{p+2i+1}=A$, then, we are still in the induction loop, since we cannot decide the value of $A$ with the equation. Indeed, we can in this case only know that $V_{p+2i+1}=V_{p-1}=V_{p+1}$. Otherwise, the equation allows us to ``collapse'' the unknown values of $V_{p-2i+1},\ldots,V_{p+2i-1}$, as well as the value of $V_{p+2i+1}$, in which case we end the induction. Recall that by induction hypothesis and after the even step, we have $V_{p-2i}=V_{p-2i+4}=\cdots=V_{p+2i}=B$ and $V_{p-2i+1}=V_{p-2i+3}=\cdots=V_{p+2i-1}=A$. Similarly, using \Cref{clm-evenBranch-1}, the induction hypothesis, and the facts that $V_{p+2i+1}=A$ and $V_{p-2i}=B$ (we have $1\leq i \leq \frac{p-3}{2}$):
	\begin{eqnarray*}
		S_{p-2i+1}(\mathcal{S})-S_{p-2i}(\mathcal{S}) & = & V_{p-2i-1} + 2V_{p-2i} + \sum_{j=p-2i+1}^{p+2i+1} \left( V_j \right) \\
		& = & V_{p-2i-1} + 2B + (2i+1)A + (2i)B.
	\end{eqnarray*}
	Again, if $V_{p-2i-1} = A$, then, we cannot decide their values with the equation and thus are still in the induction loop; if at least one of those conditions does not hold, then we can deduce their values as well as the unknown values of $V_{2p-2i+1},\ldots,V_{p+2i-1}$, in which case the induction ends.
	
	\smallskip
	Hence, either we are able to determine the unknown values (ending the induction at any point during the even or odd step), after which we complete the reconstruction using the same method as in Case~1 and Subcase~2.1 (we do not give the details here, since this leads to the exactly same computations, so we refer the reader to Case~1), or the induction continues until we reach the last possible value for $i$.
	
	If the induction does not stop early and we continue in this way until $i=\lfloor \frac{p-3}{2}\rfloor$, then, by (the already proved) induction hypothesis, for an odd $p$ we have $V_2=V_4=\cdots=V_{2p-2}=A$ and $V_3=V_5=\cdots=V_{2p-3}=B$. We then compute $S_{2p-1}(\mathcal{S})-S_{2p}(\mathcal{S})$ to obtain $V_{2p-1}$ and $S_{2p}(\mathcal{S})-S_{2p+1}(\mathcal{S})$ to obtain $V_{2p}$. As above, we check if these computations cause the collapse and allow us to finish the proof as in Case~1. If not, we have $V_2=V_4=\cdots=V_{2p}=A$ and $V_3=V_5=\cdots=V_{2p-1}=B$. Recall that we know the labels of the leaves, and let $L$ be the sum of those three labels. Hence, $S_1(\mathcal{S})-L-V_1=pA+(p-1)B$ and therefore, we can deduce the values of both $A$ and $B$.
	
	Consider then an even $p$. We have, by (the already proved) induction hypothesis, $V_3=V_5=\cdots=V_{2p-3}=A$ and $V_4=V_6=\cdots=V_{2p-4}=B$. We then compute first $S_{4}(\mathcal{S})-S_{3}(\mathcal{S})$ to obtain $V_{2}$. After this, we obtain values of $V_{2p-2}, V_{2p-1}$ and $ V_{2p}$ by computing $S_{2p-2}(\mathcal{S})-S_{2p-1}(\mathcal{S})$, $S_{2p-1}(\mathcal{S})-S_{2p}(\mathcal{S})$ and $S_{2p}(\mathcal{S})-S_{2p+1}(\mathcal{S})$, respectively. As above, we check if these computations cause the collapse and allow us to finish the proof as in Case~1. If not, we have $V_2=V_4=\cdots=V_{2p}=B$ and $V_3=V_5=\cdots=V_{2p-1}=A$. Hence, $S_1(\mathcal{S})-L-V_1=pB+(p-1)A$ and therefore, we can again deduce the values of both $A$ and $B$.
	
	\smallskip
	Hence, we know all the labels of $v_1,\ldots,v_{2p}$. Note that knowing the value of $V_1$ was primordial for this last part, which will be important later in \Cref{thm-triangleWithTail}. Finally, by computing $S_{2p+1}(\mathcal{S}) - S_{2p+2}(\mathcal{S}) = \sum_{i=1}^{2p-1} \left( V_i \right) - V_{2p+1}$, we can compute the label of $v_{2p+1}$, and thus the labels of $u_1$ and $u_2$ which are the two other leaves, completing the reconstruction.
\end{proof}

\subsection{Summary}
\label{subsec-SummaryOfReconstructing}

We considered several families of trees, which we were able to reconstruct. There are a few interesting remarks to draw from this. First, note that not all paths are reconstructable~\cite{acharya2015string}, but creating a false twin of one of the leaves of a path, that is appending a leaf to its second vertex, makes the graph always reconstructable (\Cref{thm-reconstructingS11k} for $S_{1,1,m}$). A natural question would be to further this by appending a leaf to other vertices of a path, studying the reconstructability of the subdivided star $S_{1,m,t}$. As we will see later (\Cref{fig-s123}), this family is not reconstructable in general, but it might be possible under some conditions (such as $m \leq \frac{t}{4}$).

However, the most interesting point is the diversity of algorithms for reconstruction that we found. As explained in \Cref{SecLabels}, one can reconstruct a reconstructable graph by enumerating all subgraph composition multisets of all possible labelings and comparing it to the subgraph composition multiset of the graph, but doing so can be extremely costly. We are able to be way more fine-grained, using only two sum-queries to reconstruct a star (\Cref{cor-reconstructingStars}), and four queries (one multiset and three sums) for stars subdivided at most once when the alphabet contains two symbols (\Cref{thm-reconstructingStarsSubdividedAtMostOnce}). For bistars, the alphabet size plays a role: with two symbols, three queries (two sums and one multiset) are enough (\Cref{thm-reconstructingBistarsK2}), which is also the case for a larger alphabet if the two non-leaves have different labels (Case 1 of \Cref{thm-reconstructingBistarsAnyK}). However, if those vertices have different labels, then, we may have to apply a brute-force method to find a specific composition (Case~2 of \Cref{thm-reconstructingBistarsAnyK}). Finally, for the star $S_{1,1,m}$, we use $m+1$ sum-queries, and may have different methods depending on the parity of $m$ (\Cref{thm-reconstructingS11k}).

These very different methods of reconstruction are particularly interesting, and we conjecture that the methods used above are all optimal, meaning that no reconstruction algorithm for those classes can use fewer queries. This also begs the question of the existence of graphs for which, without using a brute-force reconstruction, the highest possible number of stronger queries (that is, a linear number of multiset-queries) would still be needed. Similarly, the classes we studied either can be reconstructed using a constant and small number of queries, or seem to require a linear number of queries; are there classes that require a logarithmic, or a sublinear but non-constant number of queries?

Another question that arises from the study of bistars and stars subdivided at most once is the part played by alphabet size. Which other classes require different algorithms for reconstruction, and in particular a difference between a constant and a linear number of queries when the alphabet size $k\geq2$ passes certain thresholds?

\section{Smallest confusable graphs}\label{secSmallConfusable}

In this section, we enumerate the smallest confusable graphs and trees, and construct a general family of graphs with a dichotomy between confusable and sum-reconstructable members depending on parity of their order. We also give the smallest non-sum-reconstructable graphs.

\begin{theorem}
	\label{thm-smallestConfusableGraphs}
	The smallest confusable graphs have order~5. The smallest confusable tree has order~7.
\end{theorem}

\begin{figure}[h]
	\begin{center}
		\begin{tikzpicture}
			\node (one) at (0,0) {
				\begin{tikzpicture}
					\node[v] (0) at (0,0) {$B$};
					\node[v] (1) at (1,0.5) {$A$};
					\node[v] (2) at (0,1) {$B$};
					\node[v] (3) at (2,0.5) {$B$};
					\node[v] (4) at (3,0.5) {$A$};
					\draw (1)to(0)to(2)to(1)to(3)to(4);
				\end{tikzpicture}
			};
			
			\node (two) at (6,0) {
				\begin{tikzpicture}
					\node[v] (0) at (0,0) {$B$};
					\node[v] (1) at (1,0.5) {$B$};
					\node[v] (2) at (0,1) {$A$};
					\node[v] (3) at (2,0.5) {$A$};
					\node[v] (4) at (3,0.5) {$B$};
					\draw (1)to(0)to(2)to(1)to(3)to(4);
				\end{tikzpicture}
			};
		\end{tikzpicture}
		
		\caption{Two equicomposable 5-vertex labeled graphs.}
		\label{fig5verteqcomp}
	\end{center}
\end{figure}
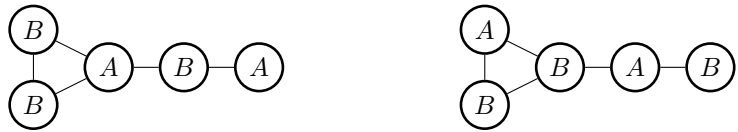

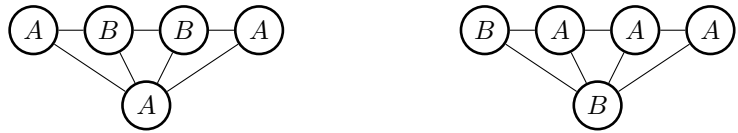
\begin{figure}[h]
	\centering
	\begin{tikzpicture}
		\node (one) at (0,0) {
			\begin{tikzpicture}
				\node[v] (0) at (0,0) {$A$};
				\node[v] (1) at (1,0) {$B$};
				\node[v] (2) at (2,0) {$B$};
				\node[v] (3) at (3,0) {$A$};
				\node[v] (4) at (1.5,-1) {$A$};
				\foreach \I in {0,1,2,3} {\draw (\I)to(4);}
				\draw (0)to(1)to(2)to(3);
			\end{tikzpicture}
		};
		
		\node (two) at (6,0) {
			\begin{tikzpicture}
				\node[v] (0) at (0,0) {$B$};
				\node[v] (1) at (1,0) {$A$};
				\node[v] (2) at (2,0) {$A$};
				\node[v] (3) at (3,0) {$A$};
				\node[v] (4) at (1.5,-1) {$B$};
				\foreach \I in {0,1,2,3} {\draw (\I)to(4);}
				\draw (0)to(1)to(2)to(3);
			\end{tikzpicture}
		};
	\end{tikzpicture}
	\caption{Two equicomposable labeled gems.}
	\label{fig-gemeqcomp}
\end{figure}

\begin{proof}
	This is an expansion on results from~\cite{acharya2015string,bartha2021reconstructibility,bartha2016reconstruction}. As in the proof of \Cref{thm-reconstructingS11k}, for sum-queries, $X$ denotes the label of vertex $x$.
	
	\medskip\noindent\textbf{Smallest confusable graphs.} No graphs of order at most~3 can be confusable, since those graphs the are paths or complete graphs.
	
	There are six connected graphs of order~4, three of them (path $P_4$, star $S_3$ and clique $K_4$) are clearly reconstructable. The cycle $C_4$ can easily be reconstructed: let $A$ be a label of a vertex, we can immediately finish if either of $A^i \in \mathcal{M}_i(\mathcal{G})$ holds for $i \in \{2,3,4\}$; otherwise, if $\mathcal{M}_4(\mathcal{G})=\{A^2XY\}$ with $X,Y \neq A$ we have two non-adjacent vertices with label $A$ and we are done, and if $\mathcal{M}_4(\mathcal{G})=\{AXYZ\}$ with $X,Y,Z \neq A$ we use $\mathcal{M}_2(\mathcal{G})$ to find the labels of the two vertices adjacent to the vertex with label $A$ and conclude with the last label. Note however that $C_4$ is not sum-reconstructable since every vertex appears with the same multiplicity in every $S_i(\mathcal{G})$.
	
	The diamond graph with vertices $V=\{u,v,w,x\}$ and edges $E=\{uv,vw,wx,xu,uw\}$, that is, complete graph $K_4$ without the edge $vx$, can easily be reconstructed. Indeed, since the query $S_2(\mathcal{G})-2S_1(\mathcal{G})$ gives the labels of $u$ and $w$, after which we can know the labels of $v$ and $x$. Vertices $u$ and $w$, and $v$ and $x$, respectively, are at symmetrical positions and hence, this is enough.
	
	Finally, the paw (a triangle with a vertex adjacent to one vertex of the triangle) can be reconstructed since the degree~3 vertex is the only one in all three subgraphs of order~3 (the other ones appear in two) so its label can be be computed with $S_3(\mathcal{G})-2S_1(\mathcal{G})$. After this, the label of the degree~1 vertex can be computed with $2S_1(\mathcal{G})-S_2(\mathcal{G})$ together with the knowledge of the label of the degree 3 vertex. Finally, we can find the labels of the two degree~2 vertices which are at symmetrical positions. Hence, no graph of order~4 is confusable.
	
	The line graph of the graph in Example~2 in~\cite{bartha2016reconstruction} (depicted on Figure~1 in their paper) is confusable, as proved in \Cref{thm-triangleWithTail} and seen in Figure~\ref{fig5verteqcomp}.
	Furthermore, the gem (a path plus one universal vertex) is confusable; consider the two labelings depicted on \Cref{fig-gemeqcomp}, we have, for both of them, $\mathcal{M}_1(\mathcal{G})=\{\{ 3A,2B \}\}$ $\mathcal{M}_2(\mathcal{G})= \{\{ 2A^2,B^2,4AB \}\}$, $\mathcal{M}_3(\mathcal{G})=\{\{ A^3,4A^2B,3AB^2 \}\}$ and $\mathcal{M}_4(\mathcal{G})=\{\{ 2A^3B,3A^2B^2 \}\}$.
	
	Let us next consider the remaining 19 connected graphs of order~5. Some can be managed easily: $K_5$, $P_5$, $S_4$ and $S_{1,1,2}$ are all reconstructable by~\cite{acharya2015string}, \Cref{cor-reconstructingStars} and \Cref{thm-reconstructingS11k}. In the following, we consider the remaining graphs of order~5, see~\cite{smallGraphs} for a list of 5-vertex graphs.
	
	\begin{enumerate}  \setcounter{enumi}{4}
		
		\item $K_5$ minus an edge. For this graph, the query $S_{2}(\mathcal{G})$ counts three times the two degree~3 vertices and four times the three degree~4 vertices, which can be separated. Since they are symmetrical, this is enough for reconstruction.
		
		\item $K_4$ plus one degree~2 vertex. Let $u_1$ and $u_2$ be the degree~3 vertices, $v_1$ and $v_2$ be the degree~4 vertices, and $w$ be the degree~2 vertex. We have $4S_3(\mathcal{G})-6S_2(\mathcal{G})=2(U_1+U_2)+8W$, from which we can clearly deduce the label of $w$, and thus the labels of $u_1$ and $u_2$, after which we can easily know the labels of $v_1$ and $v_2$.
		
		\item An edge to which three false twins are connected. The endpoints of the edge can be found using $S_2(\mathcal{G}) - 2S_1(\mathcal{G})$, after which we can deduce the labels of the three twins.
		
		\item $K_4$ plus a leaf. For this graph, the query $S_{2}(\mathcal{G})$ counts three times the three degree~3 vertices, four times the degree~4 vertex, and the leaf once. Hence, $S_2(\mathcal{G})-3S_5(\mathcal{G})$ allows us to find the label of the leaf (which is "negative") and of the degree~4 vertex. The last three vertices are twins, so this is enough.
		
		\item The cycle $C_4$ plus a vertex of degree~3. Let $u$ be the degree~2 vertex, $v_1$ and $v_2$ be its neighbours (those are false twins), and $w_1$ and $w_2$ be the other two vertices (those are true twins). We have $S_3(\mathcal{G})=5U+6(V_1+V_2)+5(W_1+W_2)$, so we can deduce the labels of $v_1$ and $v_2$. Furthermore, $S_2(\mathcal{G})$ counts three times each vertex, except $u$ which is counted twice, hence we can know its label too. This leaves the last two labels for $w_1$ and $w_2$, completing the reconstruction.
		
		\item The wheel $W_4$ (a cycle plus one universal vertex) can be reconstructed from $S_2(\mathcal{G})$: the four vertices of the cycle are symmetrical and counted three times each while the universal vertex is counted four times.
		
		\item The cycle $C_5$ can be reconstructed. Assume that the symbol $A$ is used in the labeling. Clearly, if $A^5$, $A^4$ or $A^3$ appear in $\mathcal{M}_5(\mathcal{G})$, $\mathcal{M}_4(\mathcal{G})$ or $\mathcal{M}_3(\mathcal{G})$, respectively, then we are done. If none holds, but $A^2$ appears in $\mathcal{M}_2(\mathcal{G})$, then two adjacent vertices are labeled with $A$ and their two neighbours are labeled with other symbols. If the number of $A$'s (found with $\mathcal{M}_5(\mathcal{G})$) is~3, then, the fifth vertex is labeled with $A$ and we are done. Otherwise, the query $\mathcal{M}_2(\mathcal{G})$ contains exactly two compositions $AX$ and $AY$, which give us the labels of the two neighbours, and the last vertex has the last label. Now, if there are no two consecutive $A$'s, then, there are two possibilities. Assume first that two vertices are labeled with $A$, they are at distance~2 from each other, we can find the label $X$ of their common neighbour with the composition $A^2X$ in query $\mathcal{M}_3(\mathcal{G})$, and the last two vertices are symmetric, so we are done. Assume next that each  vertex has a unique label. Let one of the vertices have label $A$. We can know the labels of its neighbours with compositions $AX$ and $AY$ in query $\mathcal{M}_2(\mathcal{G})$, after which the same query will contain compositions $XW_1$, $YW_2$ and $W_1W_2$, from which we can clearly deduce the last two labels. 
		
		\item The butterfly (two triangles attached together through a vertex $c$) can be reconstructed: the label of $c$ can be found with $S_2(\mathcal{G})-2S_1(\mathcal{G})$ since it is the only vertex appearing four times in $\mathcal{S}_2(\mathcal{G})$ (all the other ones appear twice), and then the multiset-query $\mathcal{M}_2(\mathcal{G})$ can be used to eliminate all the compositions with $c$ in it, leaving only two compositions corresponding to the pairs on each triangle.
		
		\item The kite (a cycle $vw_1xw_2$ with an edge $w_1w_2$, and a leaf $u$ attached to $v$) can be reconstructed. We have $S_3(\mathcal{G})-(S_2(\mathcal{G})+S_5(\mathcal{G})=V$, after which $2S_2(\mathcal{G})-S_3(\mathcal{G})=V+2(W_1+W_2)+X$ so we can find the labels of $w_1$ and $w_2$ (those are repeated twice, and are symmetric) and of $x$ (it appears an odd number of times). Finding the label of $u$ is then trivial.
		
		\item The dart (a cycle $vw_1xw_2$ with an edge $vx$ and a leaf $u$ attached to $v$) can be reconstructed. We have $S_2(\mathcal{G})+2S_5(\mathcal{G})-S_3(\mathcal{G}=X$, after which $S_3(\mathcal{G})-2S_2(\mathcal{G})=U-2V-2X$ so we can find the labels of $u$ and $v$. Finding the labels of $w_1$ and $w_2$ is trivial (both are false twins and we do not need to separate them).
		
		\item The house (a cycle $v_1w_1w_2v_2$ with a vertex $u$ adjacent to both $v_1$ and $v_2$). We have $S_2(\mathcal{G})-2S_5(\mathcal{G})=V_1+V_2$ and $2S_2(\mathcal{G})-S_3(\mathcal{G})=U+V_1+V_2$ so we know those three labels. All that is left is to decide which of the two remaining labels are $W_1$ and $W_2$, which we can deduce using $\mathcal{M}_2(\mathcal{G})$ (the only possible confusion would be if $W_1=V_2$ and $W_2=V_1$, which is not a problem since the two possible labelings would be isomorphic).
		
		\item The banner (a cycle $vw_1xw_2$ with a leaf $u$ attached to $v$) can be reconstructed. We have $S_4(\mathcal{G})-3S_5(\mathcal{G})=V$, after which $S_2(\mathcal{G})-2S_5(\mathcal{G})=V-U$ so we know the label of $u$ (whether $U=V$ or not). Furthermore, we have $2S_2(\mathcal{G})-S_3(\mathcal{G})=V+X$ from which we know the label of $x$, and we can find the labels of false twins $w_1$ and $w_2$.
		
		\item The cricket (a triangle $vw_1w_2$ with two leaves $u_1$ and $u_2$ attached to $v$) can be reconstructed. We have $S_2(\mathcal{G})+2S_5(\mathcal{G})-S_3(\mathcal{G})=W_1+W_2$ (both $w_1$ and $w_2$ are twins, so we know their labels), after which $3S_2(\mathcal{G})-S_3(\mathcal{G})=6V+3(W_1+W_2)$ so we can find the label of $v$. Finally, finding the labels of $u_1$ and $u_2$ is trivial (both are false twins).
		
		\item The bull (a triangle $v_1v_2w$ with two leaves $u_1$ and $u_2$ attached to $v_1$ and $v_2$, respectively) can be reconstructed. We have $2S_2(\mathcal{G})-S_3(\mathcal{G})=2(V_1+V_2)+W$, so we know the label of $w$ and the composition of labels of $v_1$ and $v_2$, then the query $\mathcal{M}_2(\mathcal{G})$ will allow us to find the labels of $u_1$ and $u_2$ (the only possible confusion is if $U_1=V_2$ and $U_2=V_1$, which is not a problem since the two possible labelings would be isomorphic).
		
		\item The complete bipartite graph $K_{2,3}$ can trivially be reconstructed from $S_2(\mathcal{G})$ (since all vertices in each part are equivalent to each other; the two parts can be separated and thus we are done).
	\end{enumerate}
	
	Hence, the only confusable graphs of order~5 are the gem and the graph depicted on \Cref{fig-triangleWithTail}, which we will call $T_2$.
	
	\medskip\noindent\textbf{Smallest confusable trees.} Trees of order at most~5 are not confusable per above.
	Trees of order~6 are $P_6$, $S_5$, $S_{1,1,3}$, the bistar $B_{2,2}$ (none of which are confusable, see \cite{acharya2015string}, \Cref{cor-reconstructingStars} and \Cref{thm-reconstructingBistarsAnyK,thm-reconstructingS11k}, respectively), and $S_{1,2,2}$ and $S_{1,1,1,2}$. For $S_{1,2,2}$, the proof of \Cref{thm-reconstructingStarsSubdividedAtMostOnce} works with any alphabet: the only possible confusion at the end of the proof (and identified in \Cref{rem-StarSubdividedAtMostOnce}) is if the two branches of length~2 are palindromes of each other, in which case the confusion does not matter due to isomorphisms. As for $S_{1,1,1,2}$, the proof of \Cref{thm-reconstructingStarsSubdividedAtMostOnce} clearly works too: once the label of the center vertex has been found, this gives the label of the only degree~2 vertex, after which we can remove pairs of leaves, hence having the composition of the branch of length~2, from which we can know the label of its leaf, giving the labels of the other three symmetrical leaves.
	
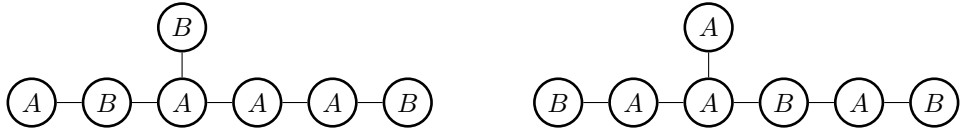
\begin{figure}[h]
	\centering
	\begin{tikzpicture}
		\node (1) at (0,0) {
			\begin{tikzpicture}
				\node[v] (0) at (0,0) {$A$};
				\node[v] (1) at (1,0) {$B$};
				\node[v] (2) at (2,0) {$A$};
				\node[v] (3) at (3,0) {$A$};
				\node[v] (4) at (4,0) {$A$};
				\node[v] (5) at (5,0) {$B$};
				\node[v] (6) at (2,1) {$B$};
				\draw (0)to(1)to(2)to(3)to(4)to(5);
				\draw (2)to(6);
			\end{tikzpicture}
		};
		
		\node (2) at (7,0) {
			\begin{tikzpicture}
				\node[v] (0) at (0,0) {$B$};
				\node[v] (1) at (1,0) {$A$};
				\node[v] (2) at (2,0) {$A$};
				\node[v] (3) at (3,0) {$B$};
				\node[v] (4) at (4,0) {$A$};
				\node[v] (5) at (5,0) {$B$};
				\node[v] (6) at (2,1) {$A$};
				\draw (0)to(1)to(2)to(3)to(4)to(5);
				\draw (2)to(6);
			\end{tikzpicture}
		};
	\end{tikzpicture}
	\caption{Two equicomposable labelings of $S_{1,2,3}$.}
	\label{fig-s123}
\end{figure}

\begin{figure}[!h]
	\centering
	\begin{subfigure}{0.45\linewidth}
		\centering
		\scalebox{0.85}{
			\begin{tikzpicture}
				\node[v] (u1) at (0,2) {$B$};
				\node[v] (u2) at (0,0) {$B$};
				\foreach \I in {1,3,5,7} {\node[v] (v\I) at (\I,1) {$A$};}
				\foreach \I in {2,4,6} {\node[v] (v\I) at (\I,1) {$B$};}
				\draw (v1)to(u1)to(u2)to(v1)to(v2)to(v3)to(v4)to(v5)to(v6)to(v7);
			\end{tikzpicture}
		}
		\caption{The labeled graph $(T_{6},\lambda_1)$.}
		\label{fig-triangleWithTail-lambda1}
	\end{subfigure}\hfill
	\begin{subfigure}{0.45\linewidth}
		\centering
		\scalebox{0.85}{
			\begin{tikzpicture}
				\node[v] (u1) at (0,2) {$A$};
				\node[v] (u2) at (0,0) {$B$};
				\foreach \I in {1,3,5,7} {\node[v] (v\I) at (\I,1) {$B$};}
				\foreach \I in {2,4,6} {\node[v] (v\I) at (\I,1) {$A$};}
				\draw (v1)to(u1)to(u2)to(v1)to(v2)to(v3)to(v4)to(v5)to(v6)to(v7);
			\end{tikzpicture}
		}
		\caption{The labeled graph $(T_{6},\lambda_2)$.}
		\label{fig-triangleWithTail-lambda2}
	\end{subfigure}
	\caption{Example of two equicomposable labelings for $T_{m}$ with even $m$.}
	\label{fig-triangleWithTail}
\end{figure}
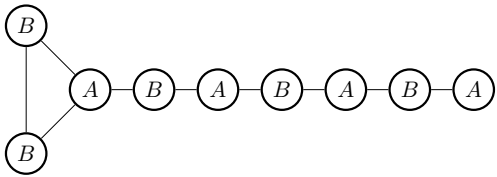
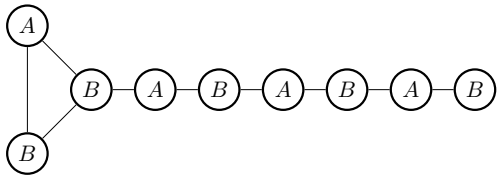

	Trees of order~7 require more work. Some of them we can easily manage: $P_7$ (\cite{acharya2015string}), the star $S_6$ (\Cref{cor-reconstructingStars}), $S_{1,1,4}$ (\Cref{thm-reconstructingS11k}) and the bistar $B_{2,3}$ (\Cref{thm-reconstructingBistarsAnyK}). The subdivided star $S_{1,2,3}$ is confusable, as can be seen on \Cref{fig-s123} (this is a result from~\cite{bartha2016reconstruction}, where they do not explicitly give the example, but Petér Burcsi and Zsuzsanna Lipták communicated the graph to us); this can be checked by enumerating the multiset compositions of both labelings. Since there are 11 trees with 7 vertices (see~\cite{smallTrees}), there are six remaining trees of order~7 to analyze.
	Note that in the following considerations, we obtain the composition of inner vertices $I$ with \Cref{lem-leavesAndInternal} and in the cases of subdivided stars (that is, in Cases~1, 2, 3 and 4) \Cref{lem-centerOfAStar} gives the label of the center vertex:
	\begin{enumerate}
		\item $S_{2,2,2}$. For this graph, the proof of \Cref{thm-reconstructingStarsSubdividedAtMostOnce} works, except in the case outlined in \Cref{rem-StarSubdividedAtMostOnce}: the branches of length~2 form a ``cycle'' (such as having compositions $AB$, $BC$ and $CA$) and we can obtain the label $X$ of the center vertex. However, in this case, the query $\mathcal{M}_4(\mathcal{S})$ will allow us to distinguish between the two possible labelings (if composition $A^2BX$ is appearing, then the branch with composition $AB$ has leaf $A$ and degree~2 vertex $B$, allowing us to find the rest; otherwise, it is the reverse).
		
		\item $S_{1,1,2,2}$. For this graph, the proof of \Cref{thm-reconstructingStarsSubdividedAtMostOnce} works: the only possible confusion at the end of the proof is if the two branches of length~2 are palindromes of each other, in which case the confusion does not matter since it leads to isomorphic labelings.
		
		\item $S_{1,1,1,3}$. Denote the center by $c$, and the vertices of the branch of length~3 by $v_1$, $v_2$ and $v_3$ where $v_1$ is adjacent to $c$ and $v_3$ is a leaf. We have $S_5(\mathcal{S})-(S_3(\mathcal{S})+S_2(\mathcal{S}))=-4C+2V_2+V_3$. Since we know $C$, this allows us to deduce $V_2$ and $V_3$, from which we can learn $V_1$ (the last internal vertex) and the three leaves adjacent to $c$.
		
		\item $S_{1,1,1,1,2}$. For this graph, there are two inner vertices. As one of them is the center vertex, we can obtain labels of both of them. Hence, we only need to identify the label of the leaf at distance two from the center.  We can do that by   removing the pairs of leaves from compositions obtained with $\bigsqcup_{X \in \mathcal{M}_5(\mathcal{S})} (S_7(\mathcal{S})-X)$, leaving only the pair with this leaf and the degree~2 vertex, whose label is known.

		\item The labeled tree $\mathcal{B}$ obtained by subdividing the middle edge of the bistar $B_{2,2}$. We can find the compositions $L$ of leaves and $I$ of internal vertices by \Cref{lem-leavesAndInternal}, from which we can compute the label of the middle vertex with $L+3I-S_2(\mathcal{B})$. Assume that the middle vertex has label $A$; we will check the largest $i$ such that the query $\mathcal{M}_i(\mathcal{B})$ contains the composition $A^i$, which will allow us to reconstruct the labeling. If $i=7$, then, each label is $A$ and we are done. If $i=6$, then, every vertex but one leaf is labeled with $A$, so we are done.
		
		If $i=5$, then, two leaves are not labeled with $A$, so we can check whether there is an $AXY$ (for $X,Y \neq A$) in $\mathcal{M}_3(\mathcal{B})$; if true, then the two non-$A$ leaves are adjacent to the same degree~3 vertex, otherwise they are adjacent to different degree~3 vertices, and in both cases we are done.
		
		If $i=4$, then, either a degree~3 vertex is labeled with $X \neq A$, or exactly three leaves are labeled with $X,Y,Z \neq A$. In the first case the two leaves adjacent to the degree~3 vertex labeled with $A$ are also labeled with $A$, and we know all other labels, so we are done. In the second case only one of $AXY$, $AXZ$ and $AYZ$ will be in $\mathcal{M}_3(\mathcal{B})$, allowing us to know which leaf is at distance~2 of the leaf labeled with $A$. In both cases we can reconstruct the labeling.
		
		If $i=3$, then, either exactly one degree~3 vertex has label $X \neq A$ and a leaf adjacent to the degree~3 vertex with label $A$ has itself label $A$ and the other leaf has label $W \neq A$, or we have $I=3A$.
		In the first case we know $X$, and we can use query $\mathcal{M}_3(\mathcal{B})$ to find all compositions $A^2Y$ with $Y \neq A$: there can be from one to five $A^2X$ depending on if zero, one or two leaves adjacent to the $X$-labeled degree~3 vertex are labeled with $A$ and whether or not $W=X$; if there are one or two $A^2X$, then $W \neq X$ and we can find $W$ from the only composition $AY$ other than $AX$ in $\mathcal{M}_2(\mathcal{B})$ after which we get the last two leaves; if there are three $A^2X$, then we have either $W=X$ and none of the last two leaves with label $A$ or $W \neq X$ and the last two leaves have label $A$, we can distinguish between them by counting the number of $AX$ compositions in $\mathcal{M}_2(\mathcal{B})$, which is two in the former case and three in the latter case, in both cases we can then easily conclude; if there are four $A^2X$, then we have $W=X$ and one of the last leaves has label $A$, we conclude immediately with the last leaf; if there are five $A^2X$, then $W=X$ and the last two leaves have label $A$. In the second case with $I=3A$, all leaves have labels different from $A$, say $W$, $X$, $Y$ and $Z$. In this case the query $\mathcal{M}_3(\mathcal{B})$ will allow us to check possible combinations among $AWX$, $AWY$, $AWZ$, $AXY$, $AXZ$ and $AYZ$; note that only two of those can exist unless there is some isomorphism, in both cases we can reconstruct the graph.
		
		If $i=2$, then, one of the degree~3 vertices has label $A$ and its two adjacent leaves have labels $X,Y \neq A$ and the other degree~3 vertex has label $Z$ (which we know). If a leaf is labeled with $A$, then it is adjacent to the $Z$-labeled degree~3 vertex. We can use query $\mathcal{M}_2(\mathcal{B})$ to get all $AX$ and $AY$ combinations, allowing us to complete the reconstruction.
		
		Finally, if $i=1$, then the two degree~3 vertices are labeled with $U,V \neq A$. Let $W$, $X$, $Y$ and $Z$ be the labels of the leaves. If $U=V$, then, we use query $\mathcal{M}_3(\mathcal{B})$. First, we count the number of $A^2U$ compositions which there are at most six; if there are six $A^2U$, then all leaves have label $A$; there cannot be exactly five $A^2U$; if there are four $A^2U$, then three leaves have label $A$ and we are done; if there are three $A^2U$, then two leaves attached to the same degree~3 vertex have label $A$; if there are two $A^2U$, then two leaves attached to different degree~3 vertices have label $A$; if there is one $A^2U$, then only one leaf has label $A$ and we may query $M_3(\mathcal{B})$ to find $WZU$ where $W,Z\neq A$ to find how labels are distributed between the leaves. In all of cases, completing the reconstruction afterwards is easy.
		If no $A^2U$ composition can be found, we remove all compositions containing $A$ and are left with the following possible compositions: $UWX$, $UWY$, $UWZ$, $UXY$, $UXZ$ and $UYZ$; only two of these exist unless there is some isomorphism, in both cases we can reconstruct the graph. If $U \neq V$, then, leaves labeled with $A$ can easily be found with compositions $AU$ and $AV$ in $\mathcal{M}_2(\mathcal{B})$ (there is at least one of each using the middle vertex, any supplementary one gives us a leaf labeled with $A$), after which we use query $\mathcal{M}_3(\mathcal{B})$ to check for $AUW$, $AVW, \ldots, AUZ$, $AVZ$, and again we can find which leaves are attached to each degree~3 vertex.
		
		\item The tree obtained by subdividing an edge incident with a leaf in the bistar $B_{2,2}$. Denote the vertices as following: $x$ is the degree~2 vertex, $y$ its adjacent leaf and $w$ its adjacent degree~3 vertex, $z$ is the leaf adjacent to $w$, $v$ is the other degree~3 vertex, and $u_1$ and $u_2$ are the two leaves adjacent to $v$. Recall that we know $L$, the composition of leaves ($u_1$, $u_2$, $y$ and $z$), and $I$, the composition of internal vertices ($v$, $w$ and $x$), thanks to \Cref{lem-leavesAndInternal}. We have $S_4(\mathcal{B})-S_3(\mathcal{B})-L-2I=Z$. Furthermore, we have $S_5(\mathcal{B})-S_3(\mathcal{B})-2L-I=V+2X$, so we can deduce those two labels, which also gives the last internal label $W$. Finally, $S_3(\mathcal{B})-2L=5V+6W+3X-Y$ allows us to know the value of $Y$, from which we finally deduce the labels of the last two leaves $u_1$ and $u_2$, completing the reconstruction.
	\end{enumerate}
	
	Hence, $S_{1,2,3}$ is the only confusable tree of order at most~7.
\end{proof}

Interestingly, one of the smallest confusable graph (illustrated in Figure \ref{fig5verteqcomp}) can be generalized into an infinite family of confusable graphs. We denote by $T_m$ the graph constructed by attaching a path of length $m$ to one vertex of a triangle (so the smallest confusable graph described above is $T_2$). We have the following:

\begin{theorem}
	\label{thm-triangleWithTail}
	Let $m$ be a positive integer and $\mathcal{T}=(T_m,\lambda)$. If $m$ is odd, then $\mathcal{T}$ is sum-reconstructable. If $m$ is even, then $T_m$ is confusable.
\end{theorem}

\begin{proof}
	Let $m$ be a positive integer and $\mathcal{T}=(T_m,\lambda)$. Denote the two degree~2 vertices in the triangle by $u_1$ and $u_2$ (note that they are twins, so if we find all the other labels then we can easily finish the reconstruction), the degree~3 vertex by $v_1$, and the vertices of the path by $v_2,\ldots,v_{m+1}$ starting from the vertex attached to the triangle. Since $T_m$ has $m+3$ vertices, we have already shown above that the 4-vertex graph with $m=1$ is not confusable.
	
	\smallskip\noindent\textbf{Case 1: $m$ is odd, so $m=2p-1$ for $p \geq 2$.} We note that, for $j\neq 2$, we have $S_j(\mathcal{T})=S_j(\mathcal{S})$ where $S_j(\mathcal{S})$ is as in Equations~(\ref{eqSi}) and~(\ref{eqSpi}) when $\mathcal{S}$ and $\mathcal{T}$ have an equal number of vertices. We further note that $S_2(\mathcal{T})\neq S_2(\mathcal{S})$. Since the proof of Lemma~\ref{lem-centerOfAStar} used $S_2(\mathcal{T})$, we cannot obtain $V_1$ as in \Cref{thm-reconstructingS11k}. First, if $p=2$, then we can reconstruct $\mathcal{T}$ in the following way: We find $V_1$ and $V_2$ with $S_3(\mathcal{T})-S_2(\mathcal{T})=V_1+2V_2$; we find $V_3$ with $S_4(\mathcal{T})-S_3(\mathcal{T})=V_3$; label $V_4$ with $S_5(\mathcal{T})-S_4(\mathcal{T})=V_4-(V_1+V_2)$; finally, we conclude with $U_1$ and $U_2$ as explained above.
	
	Now, consider $p \geq 3$. We apply the same reasoning as in Case~1 of the proof of \Cref{thm-reconstructingS11k}. Since $S_i(\mathcal{T})$ is the same for every $i \neq 2$, Equations~(\ref{eqSi}) and~(\ref{eqSpi}) still hold. The two differences are that we do not know $V_1$ and $S_2(\mathcal{T})\neq S_2(\mathcal{S})$. Hence, we can reuse \Cref{clm-oddBranch-1,clm-oddBranch-2} to alternately compute $V_{p+i}$ (using $S_{p+i+1}(\mathcal{T})-S_{p+i}(\mathcal{T})$ for increasing values of $i \in \{1,\ldots,p-1\}$), and labels $V_{i-1}$ and $V_{i}$ (using $S_{i+1}(\mathcal{T})-S_i(\mathcal{T})$ for decreasing values of $i$ for $i \in \{p,\ldots,3\}$). We refer the reader to the aforementioned proof for details; note that the two claims do not use $V_1$ so we can still apply the same reasoning. Hence, after this, we know the labels of $v_2,\ldots,v_{2p-1}$. Now, by \Cref{lem-s2G}, $S_2(\mathcal{T})=3V_1+V_{2p}+2(U_1+U_2+\sum_{i=2}^{2p-1} \left( V_i \right))$, so we have $S_3(\mathcal{T})-S_2(\mathcal{T})=V_1+2V_2+\sum_{i=3}^{2p-2} \left( V_i \right) = V_1+Q$ where $Q$ is a known quantity, so we find the label of $v_1$. Finally, as in the aforementioned proof, $S_{2p+1}(\mathcal{T})-S_{2p}(\mathcal{T}) = V_{2p} - \sum_{i=1}^{2p-2} \left( V_i \right)$, hence we obtain the label of $v_{2p}$. Finally, we conclude with the last two labels for $u_1$ and $u_2$.

	\medskip\noindent\textbf{Case 2: $m$ is even, so $m=2p$ for $p \geq 1$.} We prove that $T_m$ is confusable by exhibiting two labelings $\lambda_1$ and $\lambda_2$ such that $\mathcal{M}((T_m,\lambda_1))=\mathcal{M}((T_m,\lambda_2))$. The two binary labelings, illustrated in \Cref{fig-triangleWithTail}, are defined as follows:
	
	\begin{itemize}
		\item $\lambda_1(v_{2i+1})=A$ for $i \in \{0,\ldots,p\}$ and $\lambda_1(u_1)=\lambda_1(u_2)=\lambda_1(v_{2i})=B$ for $i \in \{1,\ldots,p\}$ (\Cref{fig-triangleWithTail-lambda1});
		\item $\lambda_2(u_1)=\lambda_2(v_{2i})=A$ for $i \in \{1,\ldots,p\}$ and $\lambda_2(u_2)=\lambda_2(v_{2i+1})=B$ for $i \in \{0,\ldots,p\}$ (\Cref{fig-triangleWithTail-lambda2}).
	\end{itemize}
	
	We next prove that $\mathcal{M}((T_m,\lambda_1))=\mathcal{M}((T_m,\lambda_2))$. First, note that both graphs have the same numbers of symbols $A$ and $B$. Furthermore, it is easy to see that, for both labeled graphs, all subgraphs of order~2 have set $\{A,B\}$ except one which has set $\{B^2\}$ (using \Cref{lem-s2G}).
	
	Consider now the subgraphs of order $2i$ for $i \in \{2,\ldots,p+1\}$.
	For $(T_m,\lambda_1)$, one subgraph has composition $A^{i-1}B^{i+1}$ (the one containing the triangle), and all other subgraphs have composition $A^iB^i$ since they are all paths of vertices labeled alternatively with symbols $A$ and $B$. The same holds for $(T_m,\lambda_2)$, where the subgraph with composition $A^{i-1}B^{i+1}$ is the one containing vertices $u_2,v_1,\ldots,v_{2i-1}$, and all other subgraphs have composition $A^iB^i$.
	
	Finally, consider the subgraphs of order $2i+1$ for $i \in \{1,\ldots,p\}$. For $(T_m,\lambda_1)$, there are $p-i+3$ subgraphs with composition $A^iB^{i+1}$ (the $p-i$ subpaths starting from $v_{2j}$ and ending at $v_{2j+2i}$ for $j \in \{1,\ldots,p-i\}$, as well as the three subgraphs containing $u_1$ and/or $u_2$) and $p-i+1$ subgraphs with composition $A^{i+1}B^i$ (the subpaths starting from $v_{2j+1}$ and ending at $v_{2j+1+2i}$ for $j \in \{0,\ldots,p-i\}$). The same holds for $(T_m,\lambda_2)$: there are $p-i+3$ subgraphs with composition $A^iB^{i+1}$ (the $p-i+1$ subpaths starting from $v_{2j+1}$ and ending at $v_{2j+1+2i}$ for $j \in \{0,\ldots,p-i\}$, as well as the two subgraphs containing $u_2$) and $p-i+1$ subgraphs with composition $A^{i+1}B^i$ (the $p-i$ subpaths starting from $v_{2j}$ and ending at $v_{2j+2i}$ for $j \in \{1,\ldots,p-i\}$, as well as the subgraph containing $u_1$ and not $u_2$).
	
	Hence, $\mathcal{M}((T_m,\lambda_1))=\mathcal{M}((T_m,\lambda_2))$, and as such $T_m$ is confusable.
\end{proof}

Compare \Cref{thm-triangleWithTail} with \Cref{thm-reconstructingS11k}: adding one edge (alternately, creating a twin instead of a false twin) made half of the graph class confusable. This is due to the fact that, in a subdivided star, finding the label of the center vertex is easy. For $T_m$ (as with $S_{1,1,m}$), we are starting the reconstruction from the middle of the path. In the odd case, this is possible since there are three middle labels $V_{p-1}, V_p$ and $V_{p+1}$ which we can reconstruct. However, in the even case, there might be some confusion with the labels of the middle vertices (even- and odd-indexed vertices alternating labels as in \Cref{fig-triangleWithTail}); but in the case of the subdivided star, this is solved by reaching an imbalance between even- and odd-indexed vertices by knowing the label of the center of the star, while $T_m$ does not allow us to decide between the two possible labelings.

Finally, we show that while the smallest confusable trees have seven vertices, the path $P_4$ on four vertices is not sum-reconstructable.
\begin{theorem}\label{thmSmallnonsumrec}
	The smallest non-sum-reconstructable graphs and trees have four vertices. 
\end{theorem}
\begin{proof}
	There are two three vertex graphs, path $P_3$ and cycle $C_3$. For $P_3$, we may obtain the center vertex with sum-compositions of size 2 which is enough to reconstruct $P_3$. Cycle $C_3$ is trivially reconstructable using sum-compositions of size $1$.
	
	Let us next consider path $P_4$ where vertices are labeled as $V_1,V_2,V_3,V_4$ from left to right. Let us denote by $S_i$ the sum-composition of size $i$. We have $S_1=S_4=V_1+V_2+V_3+V_4$, $S_2=S_3=V_1+2V_2+2V_3+V_4$. Hence, we cannot distinguish for example between $V_1=A,V_2=A,V_3=B,V_4=B$ and $V_1=A,V_2=B,V_3=A,V_4=B$.
\end{proof}

\section{Enumerating confusable graphs}\label{secConfusable}

In \cite{acharya2015string}, Acharya \emph{et al.} have presented a way to construct equicomposable pairs of paths of length $n\geq8$ where $n=m\cdot m'-1$ for  integers $m,m'\geq3$. We use their construction to give more complex examples of equicomposable graphs. Consider two labeled paths $\mathcal{P}_1$ and $\mathcal{P}_2$ where vertices of $\mathcal{P}_2$  are denoted from left to right by $p_1,p_2,\dots, p_m$ and they are labeled as $\lambda_2(p_i)=A_i$ for each $i$. 
We call the \emph{interleaving} of a path $\mathcal{P}_1$ by bits of $\mathcal{P}_2$, denoted by $\mathcal{P}_1\circ \mathcal{P}_2$, a path where we start with a copy of path $\mathcal{P}_1$, then continue with the first label of $\mathcal{P}_2$ which is then followed by another copy of $\mathcal{P}_1$ and the second label of $\mathcal{P}_2$, and so on. In this way, we construct the labeled path with the following structure: $\mathcal{P}_1 A_1 \mathcal{P}_1 A_2 \mathcal{P}_1\dots \mathcal{P}_1 A_m \mathcal{P}_1$.

We further denote by $\mathcal{P}_2^*$ the reversal of path $\mathcal{P}_2$, that is, the path labeled from left to right by $A_m,A_{m-1},\dots, A_1$. So, if $\lambda_2(p_i)=A_i$, then we denote $\lambda_2^*(p_i)=A_{m-i+1}$).

\begin{lemma}[\cite{acharya2015string}]\label{lemConfPaths}
	For any labeled paths $\mathcal{P}_1$ and $\mathcal{P}_2$, we have $\mathcal{P}_1\circ \mathcal{P}_2\sim \mathcal{P}_1\circ \mathcal{P}_2^*$.
\end{lemma}

In particular, Achariya et al. \cite{acharya2015string} proved that there exists a bijection $g$ pairing the subgraphs of $\mathcal{P}_1\circ \mathcal{P}_2$ and $\mathcal{P}_1\circ \mathcal{P}_2^*$ which have identical composition multisets. Furthermore, when $\mathcal{P}_2$ contains an odd number of vertices $m=2t-1$, we can construct the bijection $g$ so that any subgraph of $\mathcal{P}_1\circ \mathcal{P}_2$ containing the center vertex $p_t$ of $\mathcal{P}_2$ is paired with a subgraph $\mathcal{P}_1\circ \mathcal{P}_2^*$ also containing the center vertex $p_t$. We call this bijection $g$ the \emph{path bijection between $\mathcal{P}_1\circ \mathcal{P}_2$ and $\mathcal{P}_1\circ \mathcal{P}_2^*$}.

\begin{theorem}
	\label{thmConfGraphs}
	Let $\mathcal{G}$ be a labeled graph, $\mathcal{P}_1$ be a labeled path on at least two vertices, and $\mathcal{P}_2$ be a labeled path on $2t-1$ vertices $p_1,p_2,\dots,p_{2t-1}$.
	Let $\mathcal{G}_1$ be constructed by joining $\G$ with a single edge to vertex $p_t$ in $\mathcal{P}_1\circ \mathcal{P}_2$ and $\mathcal{G}_2$ be constructed by joining $\G$ with a single edge to vertex $p_t$ in $\mathcal{P}_1\circ \mathcal{P}_2^*$.
	We have $$\G_1\sim \G_2.$$
\end{theorem}

\begin{proof}
	We construct a bijection $f$ from the set of labeled induced subgraphs of $\mathcal{G}_1$ to the set of labeled induced subgraphs of $\mathcal{G}_2$ so that for any induced subgraph $\mathcal{G}'$ of $\mathcal{G}_1$ we have  $c(\mathcal{G}')=c(f(\mathcal{G}'))$. Let us denote by $x$ the vertex of $\G$ adjacent to the center vertex $p_t$ of $\Pc_2$.
	
	When subgraph $\mathcal{G}'$ contains only vertices from $\mathcal{G}$, we let $f(\mathcal{G}')=\mathcal{G}'$. When $\mathcal{G}'$ contains vertices only from $\mathcal{P}_1\circ \mathcal{P}_2$, we have $f(\mathcal{G}')=g(\mathcal{G}')$ where $g$ is the path bijection between $\mathcal{P}_1\circ \mathcal{P}_2$ and $\mathcal{P}_1\circ \mathcal{P}_2^*$. 
	
	Consider finally the case where $\G'$ contains vertices from both $\G$ and $\Pc_1\circ\Pc_2$. We partition the subgraph $\mathcal{G}'$ into two parts $\mathcal{G}'_g$ and $\mathcal{G}'_p$ where $\mathcal{G}'_g$ is a subgraph of $\mathcal{G}$ and $\mathcal{G}'_p$ is a subgraph of $\mathcal{P}_1\circ \mathcal{P}_2$. We construct  $f(\mathcal{G}')$ by separately having $f(\mathcal{G}'_g)=\mathcal{G}'_g$ and $f(\mathcal{G}'_p)=g(\mathcal{G}'_p)$ and after this, by adding an edge  between $f(\mathcal{G}'_g)$ and $f(\mathcal{G}'_p)$ from $p_t$ to $x$ (note that this is always possible since, by construction, $x \in f(\mathcal{G}'_g)$ and $p_t \in f(\mathcal{G}'_g)$).
	
	It is clear that this construction guarantees $c(f(\G'))=c(\G')$.
	Since this holds for any connected induced subgraph $\G'$ of $\G_1$, and furthermore $f(\G')$ is a  connected induced subgraph of $\G_2$, we have $\M(\G_1)=\M(\G_2)$.
\end{proof}

In the following theorem we show that while some subdivided stars are reconstructable, it is also possible find subdivided stars which contain a large number of equicomposable labelings. 

\begin{theorem}
	\label{thmEquicomposableSubdividedStars}
	For $k\geq2$ and $m\geq3$, there exist subdivided stars on $n=(2m+1)(k^m-k^{\lceil m/2\rceil})+1$ vertices which have $2^{\Theta(n/\log n)}$ equicomposable labeled graphs.
\end{theorem}

\begin{proof}
	Let $P^m$ be the set of non-isomorphic and non-palindromic labeled paths of length $m\geq2$ over an alphabet of size $k$. We note that by \Cref{thmPathLabels}, there are $\rchi_k(P_m)=\frac{k^m+k^{\lceil m/2\rceil}}{2}$ non-isomorphic labelings. Furthermore, as we have seen in the proof of \Cref{thmPathLabels}, there are $k^{\lceil m/2\rceil}$ palindromes of length $m$. Hence, we have $|P^m|=\frac{k^m-k^{\lceil m/2\rceil}}{2}$. Let $\Pc$ be a labeled path on three vertices labeled as $A$, $A$ and $B$ from left to right. For each $\Pc_m\in P^m$, we construct interleaved path $\Pc i=\Pc_m\circ \Pc$. After this, we identify the center vertex of each interleaved path $\mathcal{P}i$ leading to a subdivided star $\mathcal{S}$. Note that $\Pc i$ has $4m+3$ vertices and hence $\mathcal{S}$ has $(4m+2)|P^m|+1=n$ vertices. 
	
	Observe that, by \Cref{thmConfGraphs}, $\mathcal{S}$ is equicomposable with any subdivided star in which some $\mathcal{P}i=\mathcal{P}_m\circ \mathcal{P}$ is replaced by $\mathcal{P}_m\circ \mathcal{P}^*$. Note that since paths in $P^m$ are chosen to be non-isomorphic and non-palindromic, the obtained subdivided star is non-isomorphic with $\mathcal{S}$. Furthermore, we may replace any number of $|P^m|$ paths in this way. In other words, we have $|P^m|$ binary choices and thus, there is a class of equicomposable graphs containing $2^{|P^m|}$ labeled graphs. Since $n=(4m+2)|P^m|+1$, we have $2^{|P^m|}\in 2^{\Theta(n/\log n)}$.
\end{proof}

\section{Conclusion}

In this article, we have expanded the problem of reconstructing vertex-labeled graphs from their composition multisets from paths in~\cite{acharya2015string} to general graphs (which has been briefly considered in~\cite{bartha2016reconstruction}).  

In particular, we have given multiple reconstructable graph classes and shown reconstruction algorithms using as few queries as possible. On the negative side, we also found some non-reconstructable graph classes, and studied smallest confusable graphs. While our focus has been on trees (and particularly on subdivided stars), we also studied other classes such as paths where one leaf has been twinned.

In general, the problem of deciding whether a graph class is reconstructable seems hard, even when we want the graph classes to be rather simple for possible applications in polymer-based information storage systems. From this perspective, the graph class of subdivided stars $S_{1,1,k}$ seems quite enticing, it is at the same time simple to construct, easy to reconstruct due to sum-reconstructability (by \Cref{thm-reconstructingS11k}) and allows for a larger number of labelings than simple paths (by \Cref{thmPathLabels,thmSubStarLabels}). As an interesting open question, we ask: is any other subdivided star $S_{1,k,\ell}$ with $k>1$ (sum)-reconstructable when $\ell$ can be large? We know that the answer is negative in some cases when $k=\ell$ (by \Cref{thmConfGraphs}), so it might be worth to restrict the value of $k$ compared to $\ell$. Another possible direction for future research would be finding some conditions for graphs to be confusable. We also refer the reader back to \Cref{subsec-SummaryOfReconstructing} for discussion specifically centered on reconstructable graph classes.

\section*{Acknowledgement}

We thank Dénes Bartha, Péter Burcsi, and Zsuzsanna Lipták for personally communicating us the seven vertex confusable tree.

\bibliographystyle{plain}
\bibliography{ReconstructionBiblio}

\end{document}